\documentclass{lmcs}
\pdfoutput=1

\usepackage{lastpage}
\lmcsdoi{15}{3}{11}
\lmcsheading{}{\pageref{LastPage}}{}{}%
{Feb.~15,~2018}{Aug.~06,~2019}{}

\usepackage{amssymb,amsmath,amsfonts}
\usepackage{graphicx}
\usepackage{subfig}
\usepackage[reftex]{theoremref}
\usepackage{calc}
\newcommand{\function}[2]{#1 \to #2}

\let\pf\partialFunction
\newcommand{\mto}{\rightrightarrows}
\newcommand{\mpf}[2]{{\subseteq #1} \mto #2}
\newcommand{\mf}[2]{#1
\mto #2}

\newcommand{\dom}[1]{\mathop{\mathrm{dom}(#1)}}
\newcommand{\restr}[2]{#1 {\upharpoonright} #2}
\newcommand{\corestr}[2]{{#1} {\downharpoonright}
{#2}}
\newcommand{\preimage}[2]{#1^{-1}#2}
\newcommand{\intset}[2]{\{ #1
\,;\, #2 \}}
\newcommand{\extset}[1]{\{#1\}}
\newcommand{\tuple}[1]{\ulcorner {#1} \urcorner}
\newcommand{\nat}{\mathbb{N}}
\newcommand{\Baire}{\infseq{\nat}}
\newcommand{\Cantor}{\infseq{2}}
\newcommand{\baire}{\Baire}

\newcommand{\finbaire}{\finseq{\nat}}
\newcommand{\fininfbaire}{\fininfseq{\nat}}
\newcommand{\finseq}[1]{{{#1}^{< \nat}}}
\newcommand{\infseq}[1]{{{#1}^{\nat}}}
\newcommand{\fininfseq}[1]{#1^{\leq \nat}}
\newcommand\I{\mathrm{I}}
\newcommand\II{\mathrm{II}}
\newcommand{\repSpace}[1]{\mathbb{#1}}
\newcommand{\representation}[1]{\delta_\repSpace{#1}}
\newcommand{\transparentcylinder}[1]{#1^{\mathrm{tc}}}

\newcommand{\Wadge}{\mathcal{W}}
\newcommand{\Weihrauch}{\mathrm{W}}
\newcommand{\lseq}{\langle}
\newcommand{\rseq}{\rangle}
\newcommand{\seq}[1]{{\lseq #1 \rseq}}
\newcommand{\emptyseq}{\seq{}}
\newcommand{\intseq}[2]{\lseq #1 \,;\, #2 \rseq}
\newcommand{\last}[1]{{\bot}({#1})}
\newcommand{\reals}{\mathbb{R}}
\newcommand{\argument}{\,\cdot\,}
\newcommand{\ordinalup}[1]{{#1}\scalebox{0.7}{$\uparrow$}
}
\newcommand{\twoSpace}{\underline{2}}
\newcommand{\UC}{\mathrm{UC}}
\newcommand{\para}[1]{\widehat{#1}}
\newcommand{\WitnessleqName}{\mathrm{Neg}}
\newcommand{\Witnessleq}[1]{\WitnessleqName(#1)}
\newcommand{\WitnessleqMark}{\mathrm{N}}
\newcommand{\WitnessgeqName}{\mathrm{Pos}}
\newcommand{\Witnessgeq}[1]{\WitnessgeqName(#1)}
\newcommand{\WitnessgeqMark}{\mathrm{P}}
\newcommand{\bij}{\mathrm{e}_\nat}

\let\bisim\bisimilar
\newcommand{\PDName}{\mathrm{PD}}
\newcommand{\PD}[1]{\PDName(#1)}
\newcommand{\iPDName}{\PDName^\star}
\newcommand{\iPD}[2]{\iPDName(#1,#2)}
\newcommand{\cmfName}{\mathcal{M}}
\newcommand{\cmf}[2]{\cmfName(#1,#2)}
\newcommand{\leqWadge}{\leq_\Wadge}
\newcommand{\borelSigma}[1]{\boldsymbol{\Sigma}^0_{#1}}
\newcommand{\borelPi}[1]{\boldsymbol{\Pi}^0_{#1}}
\newcommand{\id}{\mathrm{id}}
\newcommand{\length}[1]{{\lvert{#1}\rvert}}
\newcommand{\concat}{\raisebox{1ex}{\scalebox{.6}{$\frown$}}}
\newcommand{\rank}{\mathrm{rk}}
\newcommand{\constantseq}[1]{#1^\nat}
\newcommand{\leqW}{\leq_{\Weihrauch}}
\newcommand{\leqsW}{\leq_{\mathrm{s}\Weihrauch}}
\newcommand{\leqcW}{\leqW^\mathrm{t}}
\newcommand{\leqcsW}{\leqsW^\mathrm{t}}
\let\leqscW\leqcsW

\newcommand{\equivW}{\mathrel{\equiv_{\Weihrauch}}}
\newcommand{\equivsW}{\mathrel{\equiv_{\mathrm{s}\Weihrauch}}}
\newcommand{\nleqcW}{\nleqW^\mathrm{t}}
\newcommand{\nleqW}{\nleq_\Weihrauch}
\newcommand{\comp}{\mathop\circ}
\newcommand{\starcomp}{\mathop\star}
\newcommand{\tightens}{\preceq}
\newcommand{\tightenedby}{\succeq}
\newcommand{\realizes}{\vdash}
\newcommand{\Label}{\operatorname{Label}}
\newcommand{\Linearize}{\operatorname{Prune}}
\newcommand{\shift}[1]{\mathrm{shift}({#1})}
\newcommand{\LT}{\repSpace{LT}}
\newcommand{\AT}{\repSpace{AT}}
\newcommand{\IT}{\repSpace{IT}}
\newcommand{\UT}{\repSpace{UT}}
\newcommand{\COrd}{\repSpace{CO}}
\newcommand{\deltaCOrd}{\delta_{\mathrm{nK}}}
\newcommand{\enumTrees}{\mathrm{e}_\LT}
\newcommand{\sizeName}{\mathrm{size}}
\newcommand{\size}[1]{\sizeName(#1)}
\newcommand{\open}[1]{\mathcal{O}(#1)}
\newcommand{\SubTrees}{\mathrm{SubTrees}}
\newcommand{\ConstructTree}{\mathrm{ConsTree}}
\newcommand{\graftName}{\mathrm{Graft}}
\newcommand{\graft}[1]{\graftName(#1)}
\newcommand{\val}{\mathrm{Val}}
\newcommand{\WitnessAbsence}{\operatorname{WitnessAbsence}}
\newcommand{\TreeWithRank}{\mathrm{TreeWithRank}}
\newcommand{\AuxName}{\mathrm{Aux}}
\newcommand{\Aux}[1]{\AuxName(#1)}
\newcommand{\auxName}{\mathrm{aux}}

\newcommand{\transparentName}{\mathrm{Trans}}
\newcommand{\transparent}[2]{\transparentName(#1,#2)}
\newcommand{\fb}{\mathrm{fb}}
\newcommand{\lin}{\mathrm{lin}}
\newcommand{\WFpart}{\mathrm{WF}}
\newcommand{\tree}{\Upsilon}
\newcommand{\isomorphic}{\simeq}
\newcommand{\finiteTreeProductNode}[1]{\tuple{#1}}
\newcommand{\finiteTreeProductName}{\bigotimes}
\newcommand{\finiteTreeProduct}[2]{\finiteTreeProductName_{#1}
#2}
\newcommand{\binaryTreeProductName}{\otimes}
\newcommand{\binaryTreeProduct}[2]{#1
\mathop{\binaryTreeProductName}
#2}
\newcommand{\infiniteTreeProductNode}[1]{[#1]}
\makeatletter
\newcommand*{\infiniteTreeProductName}{\DOTSB\mathop{\vphantom{\bigoplus}\mathpalette\mattt@op\relax}\slimits@}
\newcommand\mattt@op[2]{\vcenter{\m@th\hbox{\resizebox{\widthof{$#1\bigotimes$}}{!}{$\boxtimes$}}}}
\makeatother
\newcommand{\infiniteTreeProduct}[2]{\infiniteTreeProductName_{#1}
#2}
\newcommand{\finiteMixName}{\bigoplus}
\newcommand{\finiteMix}[2]{\finiteMixName_{#1}
#2}
\newcommand{\binaryMixName}{\oplus}
\newcommand{\binaryMix}[2]{#1
\mathop{\binaryMixName} #2} \makeatletter
\newcommand*{\infiniteMixName}{\DOTSB\mathop{\vphantom{\bigoplus}\mathpalette\matt@op\relax}\slimits@}
\newcommand\matt@op[2]{\vcenter{\m@th\hbox{\resizebox{\widthof{$#1\bigoplus$}}{!}{$\boxplus$}}}}
\makeatother
\newcommand{\infiniteMix}[2]{\infiniteMixName_{#1}
#2}
\newcommand{\concatenationSubtreeName}{\mathrm{Conc}}
\newcommand{\concatenationSubtree}[2]{\concatenationSubtreeName({#1},{#2})}
\newcommand{\extensionSetName}{\mathrm{Ext}}
\newcommand{\extensionSet}[2]{\extensionSetName({#1},{#2})}
\newcommand{\characteristicFunction}[1]{\chi_{{#1}}}

\makeatletter
\newenvironment{multiline*}{\begingroup\setlength{\jot}{0em}\start@multline\st@rredtrue}{\endmultline\endgroup}

\makeatother

\begin{document}

\title{Game characterizations and lower cones in
the Weihrauch degrees}
\author[H.~Nobrega]{Hugo Nobrega}
\address{Institute for Logic, Language, and
Computation, University of Amsterdam, The Netherlands}
\address{Departamento de Ci\^{e}ncia da
Computa\c{c}\~{a}o, Universidade Federal do Rio
de Janeiro, Brazil}
\email{hugonobrega@dcc.ufrj.br}

\author[A.~Pauly]{Arno Pauly}
\address{Department of Computer Science, Swansea
University, United Kingdom}
\email{a.m.pauly@swansea.ac.uk}

\thanks{This research was partially done whilst
the authors were visiting fellows at the Isaac
Newton Institute for Mathematical
Sciences in the programme `Mathematical,
Foundational and Computational Aspects of the
Higher Infinite'.
The research benefited from the Royal Society
International Exchange Grant \emph{Infinite games
in logic and Weihrauch degrees}.
The first author was partially supported by a
CAPES Science Without Borders grant (9625/13-5),
and the second author was partially
supported by the ERC inVEST (279499) project.
We are grateful to Benedikt L\"owe, Luca Motto
Ros, Takayuki Kihara and Rapha\"el Carroy for
helpful and inspiring discussions.
The contents of Section
\ref{sec:generalizedgames} of this paper has
already appeared in the conference proceedings
paper~\cite{nobrega_pauly_cie}.}

\begin{abstract}
	We introduce a parametrized version of the
	Wadge game for functions and show that each
	lower cone in the Weihrauch degrees is
	characterized by such a game.
	These parametrized Wadge games subsume the
	original Wadge game, the eraser and backtrack
	games as well as Semmes's tree games.
	In particular, we propose that the lower cones
	in the Weihrauch degrees are the answer to
	Andretta's question on which classes of
	functions admit game characterizations.
	We then discuss some applications of such
	parametrized Wadge games.
	Using machinery from Weihrauch reducibility
	theory, we introduce games characterizing every
	(transfinite) level of the Baire hierarchy
	via an iteration of a pruning derivative on
	countably branching trees.
\end{abstract}

\maketitle

\section{Introduction}

The interest in characterizations of classes of
functions in descriptive set theory via infinite
games began with a re-reading of the
seminal work of Wadge \cite{wadge_phd}, who
introduced a game in order to analyze a notion of
reducibility---\emph{Wadge
reducibility}---between subsets of the Baire
space.
In the variant---which by a slight abuse we call
the \emph{Wadge game}---two players, $\I$ and
$\II$, are given a partial function
$f:\pf{\baire}{\baire}$ and play with perfect
information for $\omega$ rounds.
In each run of this game, at each round player
$\I$ first picks a natural number and player
$\II$ responds by either picking a natural
number or passing, although she must pick natural
numbers at infinitely many rounds.
Thus, after $\omega$ rounds $\I$ and $\II$ build
elements $x \in \baire$ and $y\in \baire$,
respectively, and $\II$ wins the run if and
only if $x \not \in \dom{f}$ or $f(x) = y$.
It is an easy consequence of the original work of
Wadge that this game \emph{characterizes} the
continuous functions in the following sense.

\begin{thm}
	\th\label{wadge_game}
	A partial function $f:\function\baire\baire$ is
	relatively continuous
	iff
	player $\II$ has a winning strategy in the
	Wadge game for $f$.
\end{thm}

By giving player $\II$ more freedom in how she
builds her sequence $y \in \baire$, one can
obtain games characterizing larger classes of
functions.
For example, in the \emph{eraser game} (implicit
in \cite{duparc_veblenI}) characterizing the
Baire class 1 functions, player $\II$ is
allowed to erase past moves, the rule being that
she is only allowed to erase each position of her
output sequence finitely often.
In the \emph{backtrack game} (implicit in
\cite{wesep_wadgedegrees}) characterizing the
functions which preserve the class of
$\borelSigma{2}$ sets under preimages, player
$\II$ is allowed to erase \emph{all} of her past
moves at any given round, the rule in
this case being that she only do this finitely
many times.

In his PhD thesis \cite{semmes_phd}, Semmes
introduced the \emph{tree game} characterizing
the class of Borel functions in the Baire space.
Player $\I$ plays as in the Wadge game and
therefore builds some $x \in \baire$ after
$\omega$ rounds, but at each round $n$ player
$\II$ now plays a finite \emph{labeled tree},
i.e., a pair $(T_n,\phi_n)$ consisting of a
finite tree $T_n \subseteq \finbaire$ and a
function $\phi_n: \function{T_n {\smallsetminus}
\extset{\emptyseq}}{\nat}$, where $\emptyseq$
denotes the empty sequence.
The rules are that $T_n \subseteq T_{n+1}$ and
$\phi_n \subseteq \phi_{n+1}$ must hold for each
$n$, and that the \emph{final} labeled
tree $(T,\phi) = (\bigcup_{n \in \nat}
T_n,\bigcup_{n\in\nat}\phi_n)$ must be an
infinite tree with a unique infinite branch.
Player $\II$ then wins if the sequence of labels
along this infinite branch is exactly $f(x)$.
By providing suitable extra requirements on the
structure of the final tree, Semmes was able to
obtain a game characterizing the Baire
class 2 functions, and although this is not done
explicitly in \cite{semmes_phd}, it is not
difficult to see that restrictions of the
tree game also give his \emph{multitape} and
\emph{multitape eraser} games from
\cite{semmes_multitape}, which respectively
characterize
the class of functions which preserve
$\borelSigma{3}$ under preimages and the class of
functions for which the preimage of any
$\borelSigma{2}$ set is a $\borelSigma{3}$ set.

As examples of applications of these games,
Semmes found a new proof of a theorem of Jayne
and Rogers characterizing the class of
functions which preserve $\borelSigma{2}$ under
preimages and extended this theorem to the
classes characterized by the multitape and
multitape eraser games, by performing a detailed
analysis of the corresponding game in each case.

Given the success of such game characterizations,
in \cite{andretta_SLO} Andretta raised the
question of which classes of functions
admit a characterization by a suitable game.
Significant progress towards an answer was made
by Motto Ros in \cite{mottoros_game}: Starting
from a general definition of a reduction
game, he shows how to construct new games from
existing ones in ways that mirror the typical
constructions of classes of functions
(e.g., piecewise definitions, composition,
pointwise limits).
In particular, Motto Ros's results show that all
the usual subclasses of the Borel functions
studied in descriptive set theory admit
game characterizations.

In order to study the classes of functions
characterizable by games, we will use the
language of Weihrauch reducibility theory.
This reducibility (in its modern form) was
introduced by Gherardi and Marcone
\cite{gherardi_marcone_hahnbanach} and Brattka
and
Gherardi
\cite{brattka_gherardi_omniscience,brattka_gherardi_boundedness}
based on earlier work by Weihrauch on a
reducibility between
sets of functions analogous to Wadge reducibility
\cite{weihrauch_discontinuity,weihrauch_TTE}.

We will show that game characterizations and
Weihrauch degrees correspond closely to each
other.
We can thus employ the results and techniques
developed for Weihrauch reducibility to study
function classes in descriptive set theory,
and vice versa.
In particular, we can use the algebraic structure
available for Weihrauch degrees
\cite{higuchi_pauly_degreestructure,brattka_pauly_algebraic}
to obtain game characterizations for derived
classes of functions from game
characterizations for the original classes,
similar to the constructions found by Motto Ros
\cite{mottoros_game}.
As a further feature of our work, we should point
out that our results apply to the effective
setting firsthand, and are then lifted to
the continuous setting via relativization.
They thus follow the recipe laid out by
Moschovakis, e.g., in \cite[Section
3I]{moschovakis_book}.

While the traditional scope of descriptive set
theory is restricted to Polish spaces, their
subsets, and functions between them, these
restrictions are immaterial for the approach
presented here.
Our results naturally hold for multi-valued
functions between represented spaces.
As such, this work is part of a larger
development to extend descriptive set theory to a
more general setting, cf., e.g.,
\cite{debrecht_quasipolish,pauly_debrecht_dstcategory,pequignot_wadgesecondcountable,kihara_pauly_pointdegree,pauly_descriptivetheory}.

After recalling and preparing some notions and
results related to Weihrauch reducibility in
Section \ref{sec:weihrauch}, we introduce
our parametrized version of the Wadge game in
Section \ref{sec:generalizedgames} and discuss
applications.
In Section \ref{games_for_fixed_Baire_class}, we
introduce a notion of pruning derivative for
countably-branching trees, and show how
this gives rise to games characterizing each
(transfinite) level of the Baire hierarchy.

We will freely use standard concepts and notation
from descriptive set theory, referring to
\cite{kechris_book} for an introduction.

\section{Represented spaces and Weihrauch
reducibility}
\label{sec:weihrauch}

Represented spaces and continuous or computable
maps between them form the setting for computable
analysis.
The classical reference for computable analysis
is the textbook by Weihrauch
\cite{weihrauch_book}; for a comprehensive
introduction
more in line with the style of this paper we
refer to \cite{pauly_newintroduction}.

A \emph{represented} space $\repSpace{X} = (X,
\delta_\repSpace{X})$ is given by a set $X$ and a
partial surjection $\delta_\repSpace{X}:
\pf{\Baire}{X}$.
We will always consider $\baire$ as represented
by $\id$, and $\nat$ as represented by the
function $\delta_\nat(p) = n$ iff $p =
0^n1^\nat$.
Given a represented space $\repSpace{X}$ and $n
\in \nat$, $\repSpace{X}^n$ is the represented
space of $n$ -tuples represented in the
natural way since $\baire$ inherits particularly
nice tupling functions $\tuple{\argument}$ of all
finite arities from $\nat$.
The \emph{coproduct} of a family of represented
spaces
$\intset{\repSpace{Y}_x}{x\in\repSpace{X}}$
indexed by $\repSpace{X}$ is the
represented space $\coprod_{x\in\repSpace{X}}
\repSpace{Y}_x$ composed of pairs $(x,y)$ such
that $y \in \repSpace{Y}_x$, with the
representation given in the natural way, letting
a name for $(x,y)$ be a $\baire$ -pair of a name
for $x$ and one for $y$.
We denote by $\finseq{\repSpace{X}}$ the
represented space $\coprod_{n\in\nat}
\repSpace{X}^n$; thus $\finbaire$ can intuitively
be
seen as a represented space such that $\sigma$ is
named by $p$ iff $p$ encodes the length
$\length\sigma$ of $\sigma$ as well as the
$\length\sigma$ elements of $\sigma$ Finally,
$\infseq{\repSpace{X}}$ is the represented space
in which tuples $\seq{x_n}_{n\in\nat}$
are named by infinite tuples composed of a name
for each $x_n$ ---recall that $\baire$ has a
countable tupling function
$\tuple{\argument}:\function{\infseq{(\baire)}}{\baire}$
given by $\tuple{p_n}_{n\in\nat} = p$ iff
$p(\tuple{n,k}) = p_n(k)$ for each
$n,k \in \nat$.
This tupling function is naturally associated to
$\omega$-many corresponding projections: for each
$n \in \nat$ and $p \in \baire$, we
define $(p)_n \in \baire$ by $(p)_n(k) =
p(\tuple{n,k})$.

A (multi-valued) function between represented
spaces is just a (multi-valued) function on the
underlying sets.
We say that a partial function $F:
\pf{\Baire}{\Baire}$ is a \emph{realizer} for a
multi-valued function $f:
\mpf{\repSpace{X}}{\repSpace{Y}}$, denoted by $F
\vdash f$, if $\delta_\repSpace{Y}(F(p)) \in
f(\delta_\repSpace{X}(p))$ for all $p \in
\dom{f\delta_\repSpace{X}}$.
Then, given a class $\Lambda$ of partial
functions in $\baire$, we say
$f:\mpf{\repSpace{X}}{\repSpace{Y}}$ is in
$(\delta_\repSpace{X},\delta_\repSpace{Y})$-$\Lambda$
if it has a realizer in $\Lambda$.
When $\delta_\repSpace{X}$ and $\repSpace{Y}$ are
clear from the context, we will just say $f$ is
in $\Lambda$; thus we have
computable, continuous, etc., functions between
represented spaces.

Represented spaces and continuous functions (in
the sense just defined) generalize Polish spaces
and continuous functions (in the usual
sense).
Indeed, let $(X,\tau)$ be some Polish space, and
fix a countable dense sequence $\intseq{a_i}{i
\in \nat}$ and a compatible metric $d$.
Now define $\delta_\repSpace{X}$ by
$\delta_\repSpace{X}(p) = x$ iff $d(a_{p(i)},x) <
2^{-i}$ holds for all $i \in \nat$.
In other words, we represent a point by a
sequence of basic points converging to it with a
prescribed speed.
It is a foundational result in computable
analysis that the notion of continuity for the
represented space $(X,\delta_\repSpace{X})$
coincides with that for the Polish space
$(X,\tau)$.

\begin{defi}
	\th\label{def:weihrauch}
	Let $f$ and $g$ be partial multi-valued
	functions between represented spaces.
	We say that $f$ is \emph{Weihrauch-reducible}
	to $g$, in symbols $f\leqW g$, if there are
	computable
	functions $K:\pf{\Baire}{\Baire}$ and
	$H:\pf{\Baire\times\Baire}{\Baire}$ such that
	whenever $G \vdash g$, the function $F:= \left
	(p
	\mapsto H(p,G(K(p)))\right )$ is a realizer for
	$f$.
	If there are computable functions
	$K,H:\pf{\Baire}{\Baire}$ such that whenever $G
	\vdash g$ then $HGK \vdash f$, then we say that
	$f$
	is \emph{strongly Weihrauch-reducible} to $g$
	($f \leqsW f$).
	We write $f \leqcW g$ and $f \leqscW g$ for the
	variations where ``computable'' is replaced
	with ``continuous''.
\end{defi}

In this paper, we almost always denote (multi- or
single-valued) function composition by
juxtaposition (e.g., as we did for $HGK$ in
\th\ref{def:weihrauch}).
However, because some of the functions we use
have English words for names, for the sake of
clarity, when talking about compositions
involving such functions we will use the symbol
${\comp}$ to denote composition.

We refer the reader to
\cite{brattka_gherardi_pauly_handbook} for a
recent comprehensive survey on Weihrauch
reducibility.

A multi-valued function $f$ \emph{tightens} $g$
or is a \emph{tightening} of $g$, denoted by $f
\tightens g$, if $\dom{g} \subseteq
\dom{f}$ and $f(x) \subseteq g(x)$ whenever $x
\in \dom{g}$, cf.~\cite[Definition
7]{weihrauch_closed_under_programming}.
The following result illustrates how the notion
of tightening is a good tool for expressing
concepts in Weihrauch reducibility theory.
First, for any represented space $\repSpace{X}$,
let
$\Delta_{\repSpace{X}}:\function{\repSpace{X}}{\repSpace{X}\times\repSpace{X}}$
be
the total computable function given by
$\Delta_{\repSpace{X}}(x) = (x,x)$

\begin{prop}[Folklore; cf., e.g., {\cite[Chapter
	4]{pauly_phd}}]
	\th\label{Weihrauch_tightening}
	Let $f:\mpf{\repSpace{X}}{\repSpace{Y}}$ and
	${g:\mpf{\repSpace{Z}}{\repSpace{W}}}$.
	\begin{enumerate}
		\item
			\label{strong_Weihrauch_tightening_equality}
			The following are equivalent:
			\begin{enumerate}
				\item
					$f \leqsW g$ ($f \leqscW g$)

				\item
					there exist computable (continuous)
					$k:\mpf{\repSpace{X}}{\repSpace{Z}}$
					and
					$h:\mpf{\repSpace{W}}{\repSpace{Y}}$
					such that $f
					\tightenedby hgk$.
			\end{enumerate}

		\item
			\label{Weihrauch_tightening_equality}
			The following are equivalent:
			\begin{enumerate}
				\item
					\label{Weihrauch_tightening_equality_leqW}
					$f \leqW g$ ($f \leqcW g$);

				\item
					\label{Weihrauch_tightening_equality_tightening}
					there exist computable (continuous)
					$k:\mpf{\repSpace{X}}{\repSpace{Z}}$
					and $h:\mpf{\repSpace{X} \times
					\repSpace{W}}{\repSpace{Y}}$ such that
					$f \tightenedby h (\id_{\repSpace{X}}
					\times gk )\Delta_{\repSpace{X}}$;

				\item
					\label{Weihrauch_tightening_equality_equality}
					there exist computable (continuous)
					$k:\mpf{\repSpace{X}}{\repSpace{Z}}$
					and $h:\mpf{\repSpace{X} \times
					\repSpace{W}}{\repSpace{Y}}$ such that
					$f = h (\id_{\repSpace{X}} \times gk
					)\Delta_{\repSpace{X}}$.
			\end{enumerate}
	\end{enumerate}
\end{prop}

Although the Weihrauch degrees form a very rich
algebraic structure (cf., e.g.,
\cite{brattka_gherardi_hoelzl_probabilistic,brattka_pauly_algebraic,higuchi_pauly_degreestructure}
for surveys covering this aspect of
the Weihrauch lattice), we only need two
operations on the Weihrauch degrees,
parallelization and sequential composition.
Given a map $f:\mpf{\repSpace{X}}{\repSpace{Y}}$
between represented spaces, its
\emph{parallelization} is the map
$\para{f}:\mpf{\infseq{\repSpace{X}}}{\infseq{\repSpace{Y}}}$
given by $\seq{y_n}_{n\in\nat} \in
\para{f}(\seq{x_n}_{n\in\nat})$ iff
$y_n \in f(x_n)$ for each $n \in \nat$.
We say that $f$ is \emph{parallelizable} if $f
\equivW \para{f}$.
It is not hard to see that parallelization is a
closure operator in the Weihrauch degrees.
Rather than defining the sequential composition
operator $\starcomp$ explicitly as in
\cite{brattka_pauly_algebraic}, we will make use
of the following characterization:

\begin{thm}[Brattka \&
	Pauly~\cite{brattka_pauly_algebraic}]
	$\displaystyle f \star g \equivW \max_{\leqW}
	\intset{f' g'}{f' \leqW f \ \wedge \ g' \leqW
	g}$.
\end{thm}

\subsection{Transparent cylinders}

In this section we study properties of
\emph{transparent cylinders}, which will play a
central role in our parametrized version of the
Wadge game.

\begin{defi}[Brattka \& Gherardi
	{\cite{brattka_gherardi_omniscience}}]
	Let $f: \mpf{\repSpace{X}}{\repSpace{Y}}$.
	We call $f$
	\begin{enumerate}
		\item
			a \emph{cylinder} if $\id_\baire \times f
			\leqsW f$;
		\item
			\emph{transparent} iff for any computable
			or continuous $g:
			\mpf{\repSpace{Y}}{\repSpace{Y}}$ there is
			a computable or
			continuous, respectively, $f_g:
			\mpf{\repSpace{X}}{\repSpace{X}}$ such that
			$f f_g \tightens g f$.
	\end{enumerate}
\end{defi}

The transparent (single-valued) functions on the
Baire space were studied by de Brecht under the
name \emph{jump operator} in
\cite{debrecht_jump}.
One of the reasons for their relevance is that
they induce endofunctors on the category of
represented spaces, which in turn can
characterize function classes in DST
\cite{pauly_debrecht_synthetic}.
The term \emph{transparent} was coined in
\cite{brattka_gherardi_marcone_wkl}.
Our extension of the concept to multi-valued
functions between represented spaces is rather
straightforward, but requires the use of the
notion of tightening.
Note that if $f:
\mpf{\repSpace{X}}{\repSpace{Y}}$ is transparent,
then for every $y \in \repSpace{Y}$ there is some
$x \in \dom{f}$
with $f(x) = \extset{y}$, i.e., $f$ ~is
\emph{slim} in the terminology of
\cite[Definition
2.7]{brattka_gherardi_marcone_wkl}.

Two examples of transparent cylinders which will
be relevant in what follows are the functions
$\lim$ and $\lim_\Delta:
\pf{\baire}{\baire}$ defined by letting $\lim(p)
= \lim_{n \in \nat} (p)_n$ and letting
$\lim_\Delta(p)$ be the restriction of $\lim$ to
the domain $\intset{p \in \baire}{\exists n \in
\nat \forall k \geq n ((p)_k = (p)_n)}$.
To see a further example, related to Semmes's
tree game characterizing the Borel functions, one
first needs to define the appropriate
represented space of labeled trees.
For this, it is best to work in a quotient space
of labeled trees under bisimilarity.
The resulting quotient space can be thought of as
the space of labeled trees in which the order of
the subtrees rooted at the children
of a node, and possible repetitions among these
subtrees, are abstracted away.
Then the map $\mathrm{Prune}$, which removes from
(any representative of the equivalence class of)
a labeled tree which has one infinite
branch all of the nodes which are not part of
that infinite branch, is a transparent cylinder.
This idea will be developed in full in
Section~\ref{games_for_fixed_Baire_class} below.

\begin{prop}%
	\th\label{prop:transcyl}
	Let $f:\mpf{\repSpace{X}}{\repSpace{Y}}$ and
	$g: \mpf{\repSpace{Y}}{\repSpace{Z}}$ be
	cylinders.
	If $f$ is transparent then $g f$ is a cylinder
	and $g f \equivW g \starcomp f$.
	Furthermore, if $g$ is also transparent then so
	is $g f$.
\end{prop}
\begin{proof}
	($g f$ is a cylinder)
	As $g$ is a cylinder, there are computable $h:
	\mpf{\repSpace{Z}}{\baire \times \repSpace{Z}}$
	and $k: \mpf{\baire \times
	\repSpace{Y}}{\repSpace{Y}}$ such that
	$\id_\baire \times g \tightenedby h g k$.
	Likewise, there are computable $h':
	\mpf{\repSpace{Y}}{\baire \times \repSpace{Y}}$
	and $k': \mpf{\baire \times
	\repSpace{X}}{\repSpace{X}}$ such that
	$\id_\baire \times f \tightenedby h' f k'$.
	As composition respects tightening \cite[Lemma
	2.4(1b)]{pauly_ziegler_relative}, we conclude
	that $(\id_\baire \times g) (\id_\baire \times
	f) = \id_\baire \times (g f) \tightenedby h g k
	h' f k'$.
	Note that $k h':
	\mpf{\repSpace{Y}}{\repSpace{Y}}$ is
	computable, and as $f$ is transparent, there is
	some computable $f_{k
	h'}: \mpf{\repSpace{X}}{\repSpace{X}}$ with $k
	h' f \tightenedby f f_{k h'}$.
	But then $\id_\baire \times (g f) \tightenedby
	h g k h' f k' \tightenedby h g f f_{k h'} k'$,
	thus $h$ and $f_{k h'} k'$ witness that
	$\id_\baire \times (g f) \leqsW g f$,
	i.e., $g f$ is a cylinder.

	($g f \equivW g \starcomp f$)
	The direction $g f \leqW g \starcomp f$ is
	immediate.
	Let $f'$ and $g'$ be such that $f' \leqW f$,
	$g' \leqW g$, and $g' f'$ is defined.
	We need to show that $g' f' \leqW g f$.
	As $g$ and $f$ are cylinders, we find that
	already $g' \leqsW g$ and $f' \leqsW f$.
	Let $h,k$ witness the former and $h',k'$ the
	latter.
	We conclude $h g k h' f k' \tightens g' f'$.
	As above, there then is some computable $f_{k
	h'}$ with $k h' f \tightenedby f f_{k h'}$.
	Then $h$ and $f_{k h'} k'$ witness that $g' f'
	\leqsW g f$.

	Now suppose that $g$ is also transparent.

	($g f$ is transparent)
	Let $h: \mpf{\repSpace{Z}}{\repSpace{Z}}$ be
	computable.
	By assumption that $g$ is transparent, there is
	some computable $g_h:
	\mpf{\repSpace{Y}}{\repSpace{Y}}$ such that $g
	g_h
	\tightens h g$.
	Then there is some computable $f_{h}:
	\mpf{\repSpace{X}}{\repSpace{X}}$ with $f f_{h}
	\tightens g_h f$.
	As composition respects tightening \cite[Lemma
	2.4.1.b]{pauly_ziegler_relative}, we find that
	$h g f \tightenedby g g_h f \tightenedby g f
	f_{h}$, which is what we need.
\end{proof}

\begin{defi}%
	Given a function $f:\mpf{A}{B}$ and $C
	\subseteq B$, the \emph{corestriction} of $f$
	to $C$, denoted $\corestr{f}{C}$, is the
	restriction of $f$ to domain $\intset{x \in
	\dom{f}}{f(x) \subseteq C}$.
	This notion extends to functions between
	represented spaces in a natural way.
	A represented space $(X,\representation{X})$ is
	a \emph{subspace} of $(Y,\delta_Y)$, denoted
	$(X,\representation{X}) \subseteq
	(Y,\delta_Y)$, if $X \subseteq Y$ and
	$\representation{X} = \corestr{\delta_Y}{X}$.
\end{defi}

\begin{prop}%
	\th\label{corestr_to_larger_subspace}
	If $f:\mpf{\repSpace{X}}{\repSpace{Y}}$ and
	$\repSpace{Z} \subseteq \repSpace{W} \subseteq
	\repSpace{Y}$, then
	$\corestr{f}{\repSpace{Z}} \leqsW
	\corestr{f}{\repSpace{W}}$.
\end{prop}

\begin{prop}%
	\th\label{restricting_codomain_transparent}
	Any corestriction of a transparent map is
	transparent.
\end{prop}
\begin{proof}
	Let $f:\mpf{\repSpace{X}}{\repSpace{Y}}$ be
	transparent, and let $\repSpace{Z}$ be a
	subspace of $\repSpace{Y}$.
	Let $g:\mpf{\repSpace{Z}}{\repSpace{Z}}$ be
	computable.
	Then $g:\mpf{\repSpace{Y}}{\repSpace{Y}}$, and
	therefore there exists a computable
	$f_g:\mpf{\repSpace{X}}{\repSpace{X}}$ such
	that $f
	f_g \tightens g f$.
	Note that $\dom{g f} \subseteq
	\dom{(\corestr{f}{\repSpace{Z}}) f_g}$.
	Indeed, if $x \in \dom{g f}$, then $f f_g(x)
	\subseteq g f(x) \subseteq \repSpace{Z}$, so
	$f_g(x) \subseteq
	\dom{\corestr{f}{\repSpace{Z}}}$ and therefore
	$x \in \dom{(\corestr{f}{\repSpace{Z}}) f_g}$
	as desired.
	From this it immediately follows that
	$\restr{((\corestr{f}{\repSpace{Z}})
	f_g)}{\dom{g f}} = \restr{(f f_g)}{\dom{g
	f}}$, from which we conclude
	$(\corestr{f}{\repSpace{Z}}) f_g \tightens g
	(\corestr{f}{\repSpace{Z}})$.
\end{proof}

\begin{thm}
	\th\label{representative_Weihrauch_degree_Baire}
	Any multi-valued function between represented
	spaces is strongly-Weihrauch-equivalent to a
	multi-valued function on $\baire$.
\end{thm}
\begin{proof}
	Let $f: \mpf{\repSpace{X}}{\repSpace{Y}}$ be
	given.
	Define $f': \mpf{\baire}{\baire}$ by $\dom{f'}
	= \dom{f \representation{X}}$ and $q \in f'(p)$
	iff $\representation{Y}(q) \in f
	\representation{X} (p)$.
	To see that $f \equivsW f'$, first suppose $F
	\realizes f$, i.e., for any $p \in \dom{f
	\representation{X}}$ we have
	$\representation{Y} F(p) \in f
	\representation{X}(p)$.
	Then $F(p) \in f'(p)$, so $F \realizes f'$.
	Conversely, suppose for any $p \in \dom{f'} =
	\dom{f \representation{X}}$ we have $F(p) \in
	f'(p)$.
	But this happens iff $\representation{Y} F(p)
	\in f \representation{X} (p)$, i.e., $F
	\realizes f$.
\end{proof}

\begin{defi}
	The space $\cmf{\repSpace{X}}{\repSpace{Y}}$ of
	the \emph{strongly continuous functions}
	between represented spaces $\repSpace{X}$ and
	$\repSpace{Y}$ is defined by letting $p$ be a
	name for $f$ iff $p = 0^n1q$ and, letting $M$
	be the $n$\textsuperscript{th} Turing
	machine, we have
	\begin{enumerate}
		\item
			for every $r \in \dom{f\representation{X}}$
			and every $r' \in \baire$ we have that $M$
			computes a $\representation{Y}$-name for an
			element of $f\representation{X}(r)$ when
			given $\tuple{r,r'}$ as input and $q$ as
			oracle;

		\item
			for every $x \in \dom{f}$ and every $y \in
			f(x)$ there exist a
			$\representation{X}$-name $r$ for $x$ and
			an $r' \in \baire$ such
			that $M$ computes a
			$\representation{Y}$-name for $y$ when
			given $\tuple{r,r'}$ as input and $q$ as
			oracle;

		\item
			for every $x \in \repSpace{X}
			\smallsetminus \dom{f}$, every
			$\representation{X}$-name $r$ for $x$, and
			every $r' \in \baire$, we
			have that $M$ does not compute an element
			in $\dom{\representation{Y}}$ when given
			$\tuple{r,r'}$ as input and $q$ as oracle.
	\end{enumerate}
	In this case we also say that $f$ is
	\emph{strongly continuous}, and that $M$
	\emph{strongly computes} $f$ with oracle $p$.
	As expected, if the oracle $q \in \baire$ in
	the definition above is computable then $f$ is
	called \emph{strongly computable}.
\end{defi}

\begin{thm}[Brattka \& Pauly {\cite[Lemma
	13]{brattka_pauly_algebraic}}]
	\th\label{strongly_c_tightening}
	Every computable or continuous
	$f:\mpf{\repSpace{X}}{\repSpace{Y}}$ has a
	strongly computable or strongly continuous,
	respectively,
	tightening
	$g:\mpf{\repSpace{X}}{\repSpace{Y}}$.
\end{thm}
\begin{proof}
	We can assign to each Turing machine $M$ and
	oracle $q$ a function
	$g_{M,q}:\mpf{\repSpace{X}}{\repSpace{Y}}$
	given by $\dom{g_{M,q}}
	= \intset{\representation{X}(r)}{M$ produces an
	element of $\dom{\representation{Y}}$ when run
	on input $r$ with oracle $q}$ and
	$g_{M,q}(x) = \intset{\representation{Y}(q')}{$
	there exists a $\representation{X}$-name $r$
	for $x$ such that $q'$ is the output of
	$M$ when run with input $r$ and oracle $q}$.
	Now, if $M$ with oracle $q$ computes a realizer
	for $f:\mpf{\repSpace{X}}{\repSpace{Y}}$, then
	it immediately follows that $g_{M,q}
	\tightens f$.
	Finally, to see that $g_{M,q}$ is strongly
	continuous or strongly computable (in case $q$
	is computable), let $M'$ be the Turing
	machine which, on input $\tuple{r,r'}$ and with
	oracle $q$, simply runs the Turing machine $M$
	on input $r$ and oracle $q$.
	We now have that $M'$ strongly computes
	$g_{M,q}$ with oracle $q$.
\end{proof}

\begin{thm}[Brattka \& Pauly, implicit in
	{\cite[Section 3.2]{brattka_pauly_algebraic}}]
	\th\label{transparent_cylinder}
	Every multi-valued function $f$ is strongly
	Weihrauch-equivalent to some transparent
	cylinder $\transparentcylinder{f}$, which can
	furthermore be taken to have codomain $\baire$.
\end{thm}
\begin{proof}
	By
	\th\ref{representative_Weihrauch_degree_Baire},
	it is enough to prove the result for
	$f:\mpf{\baire}{\baire}$.
	Let
	$\transparentcylinder{f}:\mpf{\cmf{\baire}{\baire}\times\baire}{\baire}$
	be given by $\transparentcylinder{f}(h,x) = h
	f(x)$.
	That $\transparentcylinder{f} \leqsW f$ holds
	is of course immediate, and conversely we have
	$f \leqsW \transparentcylinder{f}$ since
	$\id_\baire$ has a computable name in
	$\cmf{\baire}{\baire}$, so the function $K(x) =
	(\id_\baire,x)$ is computable and $f =
	\transparentcylinder{f} K$.
	To see that $\transparentcylinder{f}$ is a
	cylinder, define computable
	$K:\pf{\baire\times(\cmf{\baire}{\baire}\times\baire)}{\baire}$
	and
	$H:\pf{\baire\times\baire}{\baire\times\baire}$
	by $K(p,(h,x)) = (h_p,x)$ where $h_p(y) =
	\tuple{p,h(y)}$ and $H(\tuple{p,y}) =
	(p,y)$.
	Then $H \transparentcylinder{f} K(p,(h,x)) =
	H(h_p f(x)) = H(\tuple{p,f(x)}) = (p,f(x))$, so
	$\id_\baire \times f \leqsW
	\transparentcylinder{f}$.
	Since $f \equivsW \transparentcylinder{f}$,
	this suffices.
	Finally, to see that $\transparentcylinder{f}$
	is transparent, let $g:\mpf{\baire}{\baire}$ be
	continuous or computable.
	Define
	$g':\mpf{\cmf{\baire}{\baire}\times\baire}{\cmf{\baire}{\baire}\times\baire}$
	by $g'(h,x) = (g h,x)$.
	Note that $g'$ is continuous or computable,
	respectively, since $g$ is.
	Furthermore, we have $\transparentcylinder{f}
	g'(h,x) = \transparentcylinder{f}(g h,x) = g h
	f(x) = g \transparentcylinder{f}(h,x)$, i.e.,
	$\transparentcylinder{f} g' = g
	\transparentcylinder{f}$ as desired.
\end{proof}

\begin{defi}
	We say that a represented space $\repSpace{X}$
	\emph{(strongly) encodes} $\baire$ if any
	$f:\mpf{\baire}{\baire}$ is (strongly)
	Weihrauch-equivalent to some
	$f':\mpf{\baire}{\repSpace{X}}$.
\end{defi}

Note that if $\repSpace{X}$ has a subspace which
is computably isomorphic to $\baire$, then
$\repSpace{X}$ strongly encodes $\baire$.

\begin{thm}
	\th\label{corestr_of_transparent_cylinder}
	Let $f:\mpf{\repSpace{X}}{\repSpace{Y}}$ be a
	transparent cylinder.
	If $\repSpace{Z} \subseteq \repSpace{Y}$
	(strongly) encodes $\baire$, then
	$\corestr{f}{\repSpace{Z}}$ is transparent and
	(strongly)
	Weihrauch-equivalent to $f$.
	In the strong case, $\corestr{f}{\repSpace{Z}}$
	is also a cylinder.
\end{thm}
\begin{proof}
	Note that $\corestr{f}{\repSpace{Z}} \leqsW f$
	holds for any $f$ and $\repSpace{Z}$, and if
	$f$ is transparent then so is
	$\corestr{f}{\repSpace{Z}}$.
	Now, by
	\th\ref{representative_Weihrauch_degree_Baire},
	there is some $g:\mpf{\baire}{\baire}$ which is
	strongly Weihrauch-equivalent
	to $f$.
	Therefore, by assumption, there exists
	$g':\mpf{\baire}{\repSpace{Z}}$ such that $g'$
	is (strongly) Weihrauch-equivalent to $f$.
	Since $f$ is a transparent cylinder, there
	exists a computable $h:\pf{\baire}{\baire}$
	such that $g' \tightenedby f h$.
	Hence, since the codomain of $g'$ is
	$\repSpace{Z}$, it follows that $g'
	\tightenedby (\corestr{f}{\repSpace{Z}}) h$,
	i.e., $g'
	\leqsW \corestr{f}{\repSpace{Z}}$ and therefore
	$f$ is (strongly) Weihrauch-reducible to
	$\corestr{f}{\repSpace{Z}}$.
	Finally, if $\repSpace{Z}$ strongly encodes
	$\baire$, then we have $\id_\baire \times
	\corestr{f}{\repSpace{Z}} \leqsW \id_\baire
	\times f \leqsW f \equivsW
	\corestr{f}{\repSpace{Z}}$, so
	$\corestr{f}{\repSpace{Z}}$ is a cylinder.
\end{proof}

\section{Parametrized Wadge games}
\label{sec:generalizedgames}
\subsection{The definition}

In order to define our parametrization of the
Wadge game, first we need the following notion,
which is just the dual notion to being an
admissible representation as in
\cite{schroeder_admissibility}.

\begin{defi}
	A \emph{probe} for $\repSpace{Y}$ is a
	computable partial function $\pi:
	\pf{\repSpace{Y}}{\Baire}$ such that for every
	computable or
	continuous $f: \mpf{\repSpace{Y}}{\Baire}$
	there is a computable or continuous,
	respectively, $e:
	\mpf{\repSpace{Y}}{\repSpace{Y}}$
	such that $\pi e \tightens f$.
\end{defi}

Note that a probe is always transparent, and that
the partial inverse of a computable embedding
from $\Baire$ into $\repSpace{Y}$ is
always a probe.
The following definition generalizes the
definition of a reduction game from
\cite[Subsection 3.1]{mottoros_game}, which is
recovered as
the special case in which all involved spaces are
$\baire$, the map $\pi$ is the identity on
$\baire$, and $\Xi$ is a single-valued
function.

\begin{defi}
	\th\label{generalized_Wadge_game}
	Let $\pi:\pf{\repSpace{Y}}{\Baire}$ be a probe
	and $\Xi: \mpf{\repSpace{X}}{\repSpace{Y}}$.
	The \emph{Wadge game parametrized by $\Xi$ and
	$\pi$}, in short the \emph{$(\Xi,\pi)$-Wadge
	game}, is played by two players, $\I$ and
	$\II$, who take turns in infinitely many
	rounds.
	At each round of a run of the game for a given
	function $f: \mpf{\repSpace{Z}}{\repSpace{W}}$,
	player $\I$ first plays a natural
	number and player $\II$ then either plays a
	natural number or passes, as long as she plays
	natural numbers infinitely often.
	Therefore, after $\omega$ rounds player $\I$
	builds $x \in \Baire$ and $\II$ builds $y \in
	\Baire$, and player $\II$ \emph{wins} the
	run of the game if $x \not \in
	\dom{f\delta_\repSpace{Z}}$, or $y \in
	\dom{\delta_\repSpace{W}\pi
	\Xi\delta_\repSpace{X}}$ and
	$\delta_\repSpace{W}\pi
	\Xi\delta_\repSpace{X}(y) \subseteq
	f\delta_\repSpace{Z}(x)$.
\end{defi}

Thus, the $(\Xi,\pi)$-Wadge game is like the
Wadge game but, instead of player $\I$ building
an element $x \in \dom{f}$ and player $\II$
trying to build $f(x)$, now player $\I$ builds a
name for some element $x \in \dom{f}$ and player
$\II$ tries to build a name for some
element $y \in \repSpace{Y}$ which is transformed
by $\pi\Xi$ into a name for an element in $f(x)$.
Intuitively, the idea is that the main
transformation is done by $\Xi$, but because
fixing a parametrized game entails fixing $\Xi$,
in
order for a fixed game to be able to deal with
functions between different represented spaces
there needs to be some map which will work
as an intermediary between the target space of
$\Xi$ and the target space, say $\repSpace{W}$,
of the function in question.
This role will be played by the computable map
$\representation{W}\pi$.

It is easy to see that, restricted to
single-valued functions on $\baire$, the original
Wadge game is the
$(\id_\baire,\id_\baire)$-Wadge game, the eraser
game is the $(\lim,\id_\baire)$-Wadge game, and
the backtrack game is the
$(\lim_\Delta,\id_\baire)$-Wadge game.
Semmes's tree game for the Borel functions is the
$(\mathrm{Prune},\Label)$-Wadge game, where
$\Label$ is the function extracting the
infinite running label from (any representative
of the equivalence class of) a pruned labeled
tree consisting of exactly one infinite
branch.
The details of this last example, including the
definitions of the represented spaces involved,
will be given in Section
\ref{games_for_fixed_Baire_class} below.

\begin{thm}
	\th\label{theo:main}
	Let $\Xi$, $\pi$, and $f$ be as in
	\th\ref{generalized_Wadge_game}, and
	furthermore suppose $\Xi$ is a transparent
	cylinder.
	Then player $\II$ has a (computable) winning
	strategy in the $(\Xi,\pi)$-Wadge game for $f$
	iff $f \leqcW \Xi$ ($f \leqW \Xi$).
\end{thm}
\begin{proof}
	($\Rightarrow$)
	Any (computable) strategy for player $\II$
	gives rise to a continuous (computable)
	function $k: \pf{\Baire}{\Baire}$.
	If the strategy is winning, then
	$\delta_{\repSpace{W}} \pi \Xi
	\delta_{\repSpace{X}} k \tightens f
	\delta_{\repSpace{Z}}$, which
	implies $\delta_{\repSpace{W}} \pi \Xi
	\delta_{\repSpace{X}} k
	\preimage\delta_{\repSpace{Z}} \tightens f
	\delta_{\repSpace{Z}}
	\preimage\delta_{\repSpace{Z}} = f$.
	Thus the continuous (computable) maps
	$\delta_\repSpace{W} \pi$ and
	$\delta_\repSpace{X}k\preimage{\delta_\repSpace{Z}}$
	witness
	that $f \leqcsW \Xi$ ($f \leqsW \Xi$).

	($\Leftarrow$)
	As $\Xi$ is a cylinder, if $f \leqcW \Xi$ ($f
	\leqW \Xi$), then already $f \leqcsW \Xi$ ($f
	\leqsW \Xi$).
	Thus, there are continuous (computable) $h, k$
	with $h \Xi k \tightens f$.
	As $\delta_\repSpace{W}
	\preimage{\delta_\repSpace{W}} =
	\id_\repSpace{W}$, we find that
	$\delta_\repSpace{W}
	\preimage{\delta_\repSpace{W}} h \Xi k
	\tightens f$.
	Now $\preimage{\delta_\repSpace{W}} h:
	\mpf{\repSpace{Y}}{\Baire}$ is continuous
	(computable), so by definition of a probe,
	there is some continuous (computable) $e:
	\mpf{\repSpace{Y}}{\repSpace{Y}}$ with
	$\delta_\repSpace{W} \pi e \Xi
	k \tightens f$.
	As $\Xi$ is transparent, there is some
	continuous (computable) $g$ with $e \Xi \succeq
	\Xi g$, thus $\delta_\repSpace{W}
	\pi \Xi g k \tightens f$.
	As $g k: \mpf{\repSpace{Z}}{\repSpace{X}}$ is
	continuous (computable), it has some
	(continuous) computable realizer $K:
	\pf{\Baire}{\Baire}$.
	By \th\ref{wadge_game}, player $\II$ has a
	winning strategy in the Wadge game for $K$, and
	it is easy to see that this strategy also
	wins the $(\Xi,\pi)$-Wadge game for $f$ for
	her.
\end{proof}

\begin{cor}
	Let $\Xi$ and $\Xi'$ be transparent cylinders.
	If the $(\Xi,\pi)$-Wadge game characterizes the
	class $\Lambda$ and the $(\Xi',\pi')$-Wadge
	game characterizes the class $\Lambda'$,
	then the $(\Xi' \Xi,\pi')$-Wadge game
	characterizes the class $\Lambda' \Lambda:=
	\intset{f g}{f \in \Lambda' \land
	g \in \Lambda}$.
\end{cor}
\begin{proof}
	If player $\II$ has a (computable) winning
	strategy in the $( \Xi' \Xi,\pi')$-Wadge game
	for
	$f:\mpf{\repSpace{A}}{\repSpace{B}}$, then by
	\th\ref{theo:main} we have $f \leqcW \Xi'\Xi$
	($f \leqW \Xi'\Xi$).
	Thus, by \th\ref{Weihrauch_tightening}, there
	exist continuous (computable)
	$k:\mpf{\repSpace{A}}{\repSpace{X}}$ and
	$h:\mpf{\repSpace{A} \times
	\repSpace{Z}}{\repSpace{B}}$ such that
	$h(\id_\repSpace{A} \times \Xi'\Xi
	k)\Delta_\repSpace{A} = f$.
	Now let $g' = h(\id_\repSpace{A} \times \Xi')$
	and $g = (\id_\repSpace{A} \times \Xi
	k)\Delta_\repSpace{A}$.
	Then $f = g'g$, and since $g' \leqcW \Xi'$ ($g'
	\leqW \Xi'$) and $g \leqcW \Xi$ ($g \leqW
	\Xi$), we have $g' \in \Lambda'$ and $g \in
	\Lambda$, as desired.
	Conversely, if $g = f' f$ with $f' \in
	\Lambda'$ and $f \in \Lambda$, then by
	\th\ref{theo:main} we have $f' \leqW \Xi'$ and
	$f
	\leqW \Xi$.
	Now, by \th\ref{prop:transcyl}, it follows that
	$f' f \leqW \Xi'\Xi$.
	Finally, since by \th\ref{prop:transcyl} we
	have that $\Xi'\Xi$ is a transparent cylinder,
	again by \th\ref{theo:main} it follows
	that $\II$ has a winning strategy in the
	$(\Xi'\Xi,\pi')$-Wadge game for $g$.
\end{proof}

We thus get game characterizations of many
classes of functions, including, e.g., ones not
covered by Motto Ros's constructions in
\cite{mottoros_game}.
For example, consider the function
$\mathrm{Sort}:\function{\Cantor}{\Cantor}$ given
by $\mathrm{Sort}(p) = 0^{n}1^\nat$ if $p$
contains
exactly $n$ occurrences of $0$ and
$\mathrm{Sort}(p) = 0^\nat$ otherwise.
This map was introduced by Carroy in
\cite{carroy_playingfirstbaire} (where it was
called $\mathrm{count}_0$) and studied by Neumann
and Pauly in
\cite{pauly_neumann_topologicalview}.
From the results in
\cite{pauly_neumann_topologicalview} it follows
that the class $\Lambda$ of total functions on
$\baire$ which are
Weihrauch-reducible to $\mathrm{Sort}$ is neither
the class of pointwise limits of functions in
some other class, nor the class of
$\mathbf{\Gamma}$-measurable functions for any
boldface pointclass $\mathbf{\Gamma}$ of subsets
of $\baire$ closed under countable
unions and finite intersections.
By \th\ref{transparent_cylinder}, $\mathrm{Sort}$
is Weihrauch-equivalent to some transparent
cylinder
$\transparentcylinder{\mathrm{Sort}}$ with
codomain $\baire$.
Thus, by \th\ref{theo:main}, $\Lambda$ is
characterized by the
$(\transparentcylinder{\mathrm{Sort}},\id_\baire)$-Wadge
game.

The converse of \th\ref{theo:main} is almost
true, as well:

\begin{prop}
	\th\label{prop:ifgamecharacterizes}
	If the $(\Xi,\pi)$-Wadge game characterizes a
	lower cone in the Weihrauch degrees, then it is
	the lower cone of $\pi \Xi$, and
	$\pi \Xi$ is a transparent cylinder.
\end{prop}
\begin{proof}
	Similar to the corresponding observation in
	\th\ref{theo:main}, note that whenever player
	$\II$ has a (computable) winning strategy in
	the $(\Xi,\pi)$-Wadge game for $f$, this
	induces a (strong) Weihrauch-reduction $f
	\leqcsW \pi \Xi$ ($f \leqsW \pi \Xi$).
	Conversely, by simply copying player $\I$ 's
	moves, player $\II$ wins the $(\Xi,\pi)$-Wadge
	game for $\pi \Xi$.
	This establishes the first claim.
	Now, as $\id_\Baire \times (\pi \Xi) \leqW \pi
	\Xi$, the assumption that the $(\Xi,\pi)$-Wadge
	game characterize a lower
	cone in the Weihrauch degrees implies that
	player $\II$ wins the $(\Xi,\pi)$-Wadge game
	for $\id_\Baire \times (\pi \Xi)$.
	Thus, $\id_\Baire \times (\pi \Xi) \leqsW \pi
	\Xi$ follows, and we find $\pi \Xi$ to be a
	cylinder.
	For the remaining claim that $\pi \Xi$ is
	transparent, let $g: \mpf{\Baire}{\Baire}$ be
	continuous (computable).
	Then $g \pi \Xi
	\leqcW \pi \Xi$ ($g \pi \Xi \leqW \pi \Xi$),
	hence player $\II$ has a (computable) winning
	strategy in the
	$(\Xi,\pi)$-Wadge game for $g \pi \Xi$.
	This strategy induces some continuous
	(computable) $H: \pf{\Baire}{\Baire}$ with $g
	\pi \Xi \delta_\repSpace{X}
	\succeq \pi \Xi \delta_\repSpace{X} H$.
	Thus, $\delta_\repSpace{X} H
	\preimage{\delta_\repSpace{X}}$ is the desired
	witness.
\end{proof}

\subsection{Using game characterizations}

One main advantage of having game
characterizations of some properties is realized
together with determinacy: either by choosing our
set-theoretic axioms accordingly, or by
restricting to simple cases and invoking, e.g.,
Borel determinacy, we can conclude that if the
property is false, i.e., player $\II$ has no
winning strategy, then player $\I$ has a winning
strategy.
Thus, player $\I$ 's winning strategies serve as
explicit witnesses of the failure of a property.
Applying this line of reasoning to our
parametrized Wadge games, we obtain the following
corollaries of \th\ref{theo:main}:

\begin{cor}[\textrm{ZFC}]
	Let $\Xi$ be a transparent cylinder and $\pi$ a
	probe such that $\pi \Xi$ is single-valued and
	$\dom{\pi \Xi}$ is Borel.
	Then for any $f:
	\mf{\repSpace{X}}{\repSpace{Y}}$ such that
	$\dom{\delta_\repSpace{X}}$ and $f(x)$ are
	Borel for any $x \in
	\repSpace{X}$, we find that $f \nleqcW \Xi$ iff
	player $\I$ has a winning strategy in the
	$(\Xi,\pi)$-Wadge game for $f$.
\end{cor}

\begin{cor}[\textrm{ZF} + \textrm{DC} +
	\textrm{AD}]
	Let $\Xi$ be a transparent cylinder and $\pi$ a
	probe.
	Then $f \nleqcW \Xi$ iff player $\I$ has a
	winning strategy in the $(\Xi,\pi)$-Wadge game
	for $f$.
\end{cor}

Unfortunately, as determinacy fails in a
computable setting (cf., e.g.,
\cite{cenzer_remmel_games,leroux_pauly_equilibria}),
we do not
retain the computable counterparts.
More generally, we lack a clear grasp on the
connections between winning strategies of player
$\I$ in the $(\Xi,\pi)$-Wadge game for a
function $f$ and \emph{positive} witnesses of the
fact that $f$ is not in the class characterized
by the game.
As pointed out by Carroy and Louveau in private
communication, this is true even for the original
Wadge game for functions, i.e., the
$(\id_\baire,\id_\baire)$-Wadge game.
Here we already have a notion of positive
witnesses for discontinuity, viz.\ \emph{points
of discontinuity}, and can therefore make this
discussion mathematically precise:

\begin{qu}
	Let a point of discontinuity of a function $f:
	\function{\Baire}{\Baire}$ be given as a
	sequence $(x_n)_{n \in \nat}$, a point $x \in
	\Baire$, and $\sigma \in \finbaire$ with
	$\sigma \subseteq f(x)$ such that $\forall n (
	d(x_n, x) < 2^{-n} \wedge \sigma \not\subseteq
	f(x_n))$.
	Let $\textrm{DiscPoint}$ be the multi-valued
	map that takes as input a winning strategy for
	player $\I$ in the
	$(\id_\baire,\id_\baire)$-Wadge game
	for some function $f:
	\function{\Baire}{\Baire}$, and outputs a point
	of discontinuity for that function.
	Is $\textrm{DiscPoint}$ computable?
	More generally, what is the Weihrauch degree of
	$\textrm{DiscPoint}$?
\end{qu}

We can somewhat restrict the range of potential
answers for the preceding question:

\begin{thm}
	\th\label{theo:computablewinningstrategy}
	Let player $\I$ have a computable winning
	strategy in the $(\id_\baire,\id_\baire)$-Wadge
	game for $f: \function{\Baire}{\Baire}$.
	Then $f$ has a computable point of
	discontinuity.
\end{thm}

The proof of the theorem will require some
recursion theoretic preparations.
Given $p, q \in \Baire$, let $[p\mid q] \in
\Baire$ be defined as $[p \mid q] = 0^{q(0)}(p(0)
+ 1)0^{q(1)}(p(1)+1)0^{q(2)}\ldots$, i.e.,
$[p\mid q]$ increases each number in $p$ by $1$,
and then intersperses zeros between the entries,
with the number of repetitions being
provided by $q$.
Now, given $r \in \Baire$ and some $A \subseteq
\Baire$, let $A^{+r}:= \intset{[p\mid q]}{p \in A
\wedge \forall n \in \nat (q(n) \geq
r(n))}$.

The proof of the following lemma is based on
helpful comments by Takayuki Kihara in personal
communication.

\begin{lem}
	\th\label{lem:takayuki}
	Let $F: \function{\Baire}{\Baire}$ be
	computable, $r \in \Baire$, $A, B \subseteq
	\Baire$, $B \neq \varnothing$ be such that
	$F[B^{+r}] \subseteq A$ and $A$ is
	$\Sigma^0_2$.
	Then $A$ contains a computable point.
\end{lem}
\begin{proof}
	Let $A = \bigcup_{n \in \nat} Q_n$ with
	$\Pi^0_1$ -sets $Q_n$.
	For the sake of a contradiction, assume that
	$A$ and thus all $Q_n$ contain no computable
	points.
	Pick some $p \in B$.

	As $F(0^\nat)$ is computable, we find
	$F(0^\nat) \notin Q_0$.
	As $Q_0$ is $\Pi^0_1$ and $F$ computable, there
	is some $m_0 \geq r(0)$ such that
	$F[0^{m_0}\Baire] \cap Q_0 = \varnothing$.
	Next, consider $F(0^{m_0}p(0)0^\nat)$.
	Again, this is a computable point, hence there
	is some $m_1 \geq r(1)$ such that
	$F[0^{m_0}p(0)0^{m_1}\Baire] \cap Q_1 =
	\varnothing$.
	We proceed in this manner to choose all $m_i$,
	and then define $q \in \Baire$ by $q(i) = m_i$.
	Note that $q \geq r$.
	Then $[p \mid q] \in B^{+r}$, but $F([p\mid q])
	\notin A$ by construction, hence we derive the
	desired contradiction and conclude that
	$A$ contains a computable point.
\end{proof}

\begin{proof}[Proof of
	\th\ref{theo:computablewinningstrategy}]
	Let us assume that player $\I$ has a winning
	strategy in the $(\id,\id)$-Wadge game for $f:
	\function{\Baire}{\Baire}$.
	We describe how player $\II$ can coax player
	$\I$ into playing a point of discontinuity of
	$f$.
	Player $\II$ starts passing, causing player
	$\I$ to produce longer and longer prefixes of
	some $p \in \Baire$.
	If player $\I$ ever produces a prefix $p_{\leq
	n_0}$ such that $\exists k_0 \ f[p_{\leq
	n}\Baire] \subseteq k_0\Baire$, then player
	$\II$ will play $k_0$, and then goes back to
	passing.
	If subsequently, there is some $n_1$, such that
	$\exists k_1 \ f[p_{\leq n_1}\Baire] \subseteq
	k_0k_1\Baire$, then player $\II$ plays
	$k_1$, and starts passing again, etc.
	If $f$ is continuous at $p$, then player $\II$
	will play a correct response to $f$, hence
	contradict the assumption that player $\I$
	is following a winning strategy.
	Thus, $p$ has to be a point of discontinuity of
	$f$.

	Note that if player $\II$ passes even more than
	necessary, this does not change the argument at
	all.
	Thus, we find that there is some non-empty set
	$B$ and $r \in \Baire$ such that the computable
	response function $S:
	\function{\Baire}{\Baire}$ maps $B^{+r}$ into
	the set of points of discontinuity of $f$.
	The latter is a $\Sigma^0_2$ -set, hence
	\th\ref{lem:takayuki} implies that it contains
	a computable point.
\end{proof}

A more convenient way of exploiting determinacy
of the $(\Xi,\pi)$-Wadge games could perhaps be
achieved if a more symmetric version
were found.
In this, we could hope for a dual principle $S$,
where for any $f$ either $f \leqW^c \Xi$ or $S
\leqW^c f$ holds.
More generally, we hope that a better
understanding of the $(\Xi,\pi)$-Wadge games
would lead to structural results about the
Weihrauch
lattice, similar to the results obtained by
Carroy on the strong Weihrauch reducibility
\cite{carroy_quasiorder}.

\subsection{Generalized Wadge reducibility}

As mentioned in the introduction, the Wadge game
was introduced not to characterize continuous
functions, but in order to reason about a
reducibility between sets.
Given $A, B \subseteq \Baire$, we say that $A$ is
\emph{Wadge-reducible} to $B$, in symbols $A
\leqWadge B$, if there exists a
continuous $F: \function{\Baire}{\Baire}$ such
that $\preimage{F}[B] = A$ (we use the notation
$\leq_\Wadge$ instead of the more
established $\leq_W$ in order to help avoid
confusion with Weihrauch reducibility, $\leqW$).
Equivalently, we could consider the multi-valued
total function $\frac{B}{A}: \mf{\baire}{\baire}$
defined by $\frac{B}{A}(x) = B$ if
$x \in A$ and $\frac{B}{A}(x) = (\baire
{\smallsetminus} B)$ if $x \notin A$.
It is easy to see that we have $A \leqWadge B$
iff $\frac{B}{A}$ is continuous.
It is a famous structural result due to Wadge
(using Borel determinacy) that for any Borel $A,
B \subseteq \baire$, either $A \leqWadge
B$ or $\baire {\smallsetminus} B \leqWadge A$.
In particular, the Wadge hierarchy on the Borel
sets is a strict weak order of width
$2$.\footnote{A relation $R$ is a \emph{strict
weak
order} on a set $X$ if there exists some ordinal
number $\alpha$ and a partition
$\intseq{X_\beta}{\beta<\alpha}$ of $X$ such that
$x
\mathrel{R} y$ holds iff $x \in X_\beta$, $y \in
X_\gamma$, and $\beta < \gamma$.
The \emph{width} of $R$ is the supremum of the
cardinalities of the parts in the partition.}

Both definitions generalize in a natural way to
the case where $A \subseteq \repSpace{X}$ and $B
\subseteq \repSpace{Y}$ for represented
spaces $\repSpace{X}$, $\repSpace{Y}$: $A
\leq_\Wadge' B$ iff there exists a continuous
$f:\function{\repSpace{X}}{\repSpace{Y}}$ such
that $A = \preimage{f}[B]$, and
$\frac{B}{A}:\mf{\repSpace{X}}{\repSpace{Y}}$ is
defined by letting $\frac{B}{A}(x) = B$, if $x
\in A$,
and $\frac{B}{A}(x) = \repSpace{Y} \smallsetminus
B$, otherwise.
It is easy to see that if $A \leq_\Wadge'$, then
$\frac{B}{A}$ is continuous, since if
$f:\function{\repSpace{X}}{\repSpace{Y}}$ is such
that $A = \preimage{f}[B]$, then any realizer of
$f$ also realizes $\frac{B}{A}$.
However, since not every continuous multi-valued
function has a continuous uniformization, the
converse does not hold in general.
As noted, e.g., by Hertling \cite{hertling_phd},
the relation ${\leq_\Wadge'}$ restricted to
$\repSpace{X} = \repSpace{Y} = \reals$
already introduces infinite antichains in the
resulting degree structure, and Ikegami showed
that in fact the partial order
$(\wp(\nat),{\subseteq}_\mathrm{fin})$ can be
embedded into that degree
structure~\cite[Theorem~5.1.2]{ikegami_phd}.
The generalization of $\frac{B}{A}$ was proposed
by Pequignot
\cite{pequignot_wadgesecondcountable} as an
alternative\footnote{While
Pequignot only introduces the notion for second
countable $T_0$ spaces, the extension to all
represented spaces is immediate.
Note that one needs to take into account that for
general represented spaces, the Borel sets can
show unfamiliar properties, e.g., even
singletons can fail to be Borel (cf.\ also
\cite{schroeder_selivanov_hierarchiesqcb0,schroeder_selivanov_hyperprojectiveqcb0,hoyrup8}).%
}.

It is a natural variation to replace
\emph{continuous} in the definition of Wadge
reducibility by some other class of functions
(ideally
one closed under composition).
Motto Ros has shown that for the typical
candidates of more restrictive classes of
functions, the resulting degree structures will
not
share the nice properties of the standard Wadge
degrees (they are \emph{bad})
\cite{mottoros_badreducibilities}.
Larger classes of functions as reduction
witnesses have been explored by Motto Ros,
Schlicht, and Selivanov
\cite{mottoros_schlicht_selivanov_wadgelike} in
the setting of quasi-Polish spaces---using the
generalization of the first definition of
the reduction.
Here, we explore the second generalization; thus,
we define $A \leqWadge B$ iff $\frac{B}{A}$ is
continuous.

\begin{defi}
	\th\label{def:generalizedwadge}
	Given a multi-valued function $\Xi$ and $A
	\subseteq \repSpace{X}$, $B \subseteq
	\repSpace{Y}$ for represented spaces
	$\repSpace{X}$
	and $\repSpace{Y}$, let $A \leq_\Xi B$ iff
	$\frac{B}{A} \leqcW \Xi$.
\end{defi}

\begin{obs}
	If $\Xi \star \Xi \equivW \Xi$, then $\leq_\Xi$
	is a quasiorder.
\end{obs}

The following partially generalizes \cite[Theorem
6.10]{mottoros_game}:

\begin{thm}
	Let $A \subseteq \repSpace{X}$ and $B \subseteq
	\repSpace{Y}$, let $\Xi:
	\mf{\repSpace{U}}{\repSpace{V}}$ be a
	transparent cylinder,
	and let $\pi: \pf{\repSpace{Y}}{\Baire}$ be a
	probe such that the $(\Xi,\pi)$-Wadge game for
	$\frac{B}{A}$ is determined.
	Then either $A \leq_\Xi B$ or $B \leqWadge
	\Baire {\smallsetminus} A$.
\end{thm}
\begin{proof}
	If player $\II$ has a winning strategy in the
	$( \Xi,\pi)$-Wadge game for $\frac{B}{A}$, then
	by \th\ref{theo:main}, we find that
	$\frac{B}{A} \leqcW \Xi$, hence $A \leq_\Xi B$.
	Otherwise, player $\I$ has a winning strategy
	in that game.
	This winning strategy induces a continuous
	function $s: \function{\baire}{\baire}$, such
	that if player $\II$ plays $y \in \baire$,
	then player $\I$ plays $s(y) \in \baire$.
	As $\Xi$ is a transparent cylinder and $\pi$ a
	probe, since $\id_\baire \leqW \Xi$, by
	\th\ref{theo:main} player $\II$ has a winning
	strategy in the $(\Xi,\pi)$-Wadge game for
	$\id_\baire$.
	This strategy induces a continuous function $t:
	\pf{\baire}{\baire}$ such that $\pi \Xi
	\representation{X} t
	\tightens \id_\baire$, and since $\id_\baire$
	is total and single-valued, we have that $t$ is
	total and $\pi \Xi \representation{X} t =
	\id_\baire$
	Now we consider $s t:
	\function{\baire}{\baire}$.
	If $\representation{Z}(x) \in A$, then if
	player $\II$ plays $t(x)$, player $\I$ needs to
	play some $s(t(x))$ such that
	$\representation{W}(s(t(x))) \notin B$.
	Likewise, if $\representation{Z}(x) \notin A$,
	then for player $\I$ to win, it needs to be the
	case that $\representation{W}(s(t(x)))
	\in B$.
	Thus, $s t$ is a continuous realizer of
	$\frac{B}{\baire {\smallsetminus} A}$, and $B
	\leqWadge \baire {\smallsetminus} A$
	follows.
\end{proof}

\begin{cor}[\textrm{ZF} + \textrm{DC} +
	\textrm{AD}]
	Suppose $\Xi \star \Xi \equivW \Xi$.
	Then $<_\Xi$ is strict weak order of width at
	most $2$.
\end{cor}

In \cite{mottoros_borelamenable}, in a different
formalism, Motto Ros has identified sufficient
conditions on a general reduction to
ensure that its degree structure is equivalent to
the Wadge degrees.
We leave for future work the task of determining
precisely for which $\Xi$ the degree structure of
$<_\Xi$ (restricted to subsets of
$\Baire$) is equivalent to the Wadge degrees, and
which other structure types are realizable.

\section{Games for functions of a fixed Baire
class}
\label{games_for_fixed_Baire_class}

\subsection{Spaces of trees}

An \emph{unlabeled tree}, or simply \emph{tree},
is a subset of $\finbaire$ closed under the
operation of taking initial segments.
We will typically denote trees by the letters $T,
S, U$ with or without sub- or superscripts.
Given a tree $T$ and $\sigma \in \finbaire$, we
denote
\begin{enumerate}
	\item
		$\concatenationSubtree{T}{\sigma}:=
		\intset{\tau \in \finbaire}{\sigma \concat
		\tau \in T}$

	\item
		$\extensionSet{T}{\sigma}:= \intset{\tau \in
		T}{\sigma \subseteq \tau}$
\end{enumerate}

We call a tree \emph{linear} if each of its nodes
has at most one child, \emph{finitely branching}
if each of its nodes has only
finitely many children, and \emph{pruned} if each
of its nodes has at least one child.
Given a tree $T$ and $\sigma \in T$ we define the
\emph{rank} of $\sigma$ in $T$, denoted by
$\rank_T(\sigma)$, by the recursion
$\rank_T(\sigma):= \sup
\intset{\rank_T(\tau)+1}{\sigma\subset\tau\in
T}$, if $\sigma$ is in the wellfounded part of
$T$, and
$\rank_T(\sigma):= \infty$ otherwise.
By letting $\infty > \alpha$ whenever $\alpha$ is
a countable ordinal and $\infty + n = \infty$ for
all $n \in \nat$, we get
$\rank_T(\sigma) = \sup
\intset{\rank_T(\tau)+1}{\sigma\subset\tau\in T}$
in either case.
Furthermore, if $\sigma \subset \tau \in T$ and
$\length\tau \geq \length\sigma + n$, then
$\rank_T(\sigma) \geq \rank_T(\tau)+n$.
The \emph{rank} of $T$ is $\rank_T(\emptyseq)$,
if $T \neq \varnothing$, or 0 if $T =
\varnothing$.

Given $s \in \fininfbaire$ with $\length{s} > 0$,
let the \emph{left shift} of $s$, denoted by
$\shift{s}$, be the unique $t \in
\fininfbaire$ such that $s = \seq{s(0)} \concat
t$.
Given $p \in \baire$ and $\sigma \in \finbaire$,
we say that $\sigma$ is a \emph{path through} $p$
if $p(0) \neq 0$ and recursively
$\shift{\sigma}$ is a path through
$(\shift{p})_{\sigma(0)}$ in case $\sigma \neq
\emptyseq$.
Let $\UT$ be the space of unlabeled trees
represented by the total function $\delta_{\UT}$
given by
\[
\delta_{\UT}(p):=
\intset{\sigma\in\finbaire}{\sigma \text{ is a
path through } p}.
\]

A \emph{labeled tree} is a pair $(T,\varphi)$
where $T$ is a tree, called the \emph{domain} of
$(T,\varphi)$, and $\varphi:\function{T
{\smallsetminus} \extset\emptyseq}{\nat}$ is
called the
\emph{labeling function} of $(T,\varphi)$.
We typically denote labeled trees by the letter
$\tree$, with or without sub- or superscripts.
A labeled tree $\tree$ is a \emph{subtree} of a
labeled tree $\tree'$, denoted $\tree \subseteq
\tree'$, if the domain of $\tree$ is a
subset of the domain of $\tree'$ and the labeling
function of $\tree$ is a restriction of that of
$\tree'$.
It is not hard to see that there exist a
computable enumeration
$\enumTrees:\function\nat\LT$ of all finite
labeled trees and a
computable function
$\sizeName:\function{\nat}{\nat}$ such that
$\enumTrees(m) \subset \enumTrees(n)$ implies $m
< n$, and such that
$\enumTrees(n)$ has exactly $\size{n}$ nodes.
We will in general overload notation from
unlabeled to labeled trees; whenever some such
notation is used without prior introduction,
the intended meaning will be intuitive.
For example, for $\tree=(T,\varphi)$ we will
write $\sigma \in \tree$ to mean $\sigma \in T$,
or $\rank_\tree(\sigma)$ instead of
$\rank_T(\sigma)$, etc.

If $\sigma \neq \emptyseq$ is a path through $p$,
then its \emph{label according to $p$} is
$p(\sigma(0))-1$, if $\length\sigma = 1$, or
the label of $\shift{\sigma}$ according to
$(\shift{p})_{\sigma(0)}$, otherwise.
Let $\LT$ be the space of labeled trees
represented by the total function $\delta_{\LT}$
given by $\delta_{\LT}(p) =
(\delta_{\UT}(p),\varphi)$, where
$\varphi(\sigma)$ is the label of $\sigma$
according to~ $p$.

Given labeled trees $\tree = (T,\varphi)$ and
$\tree' = (T',\varphi')$, a relation $B \subseteq
T \times T'$ is called a
\emph{bisimulation} between $\tree$ and $\tree'$
in case $\sigma \mathrel{B} \tau$ implies
$\length\sigma = \length\tau$ and:
\begin{flalign}
	&\varphi(\sigma) = \varphi'(\tau)
	\tag{label}\label{lab} \\
	&\forall \sigma' \in T \; ( \sigma \subset
	\sigma' \;\Rightarrow\; \exists \tau' \in T'
	(\tau \subset \tau' \wedge \sigma' \mathrel{B}
	\tau')) \tag{forth}\label{for} \\
	&\forall \tau' \in T' \; ( \tau \subset \tau'
	\;\Rightarrow\; \exists \sigma' \in T (\sigma
	\subset \sigma' \wedge \sigma' \mathrel{B}
	\tau')) \tag{back}\label{bac}
	\\
	&\length\sigma = \ell+1 \;\Rightarrow\;
	\restr{\sigma}{\ell} \mathrel{B}
	\restr{\tau}{\ell} \tag{parent}\label{par}
\end{flalign}
It is easily seen that the union of any family of
bisimulations between given labeled trees is also
a bisimulation between those trees.
Therefore, between any pair of labeled trees
$\tree$ and $\tree'$ there always exists a
largest bisimulation, denoted
${\bisim_{\tree,\tree'}}$.
We say $\tree$ and $\tree'$ are \emph{bisimilar},
denoted $\tree \bisim \tree'$, in case
$\bisim_{\tree,\tree'}$ is nonempty.
A particular case of a bisimulation between
$\tree=(T,\varphi)$ and $\tree'=(T',\varphi')$ is
an \emph{isomorphism} between those trees,
i.e., a bijection $\iota:\function{T}{T'}$
satisfying, for any $\sigma,\tau \in T$:
\begin{enumerate}
	\item $\sigma\subseteq\tau$ iff
		$\iota(\sigma)\subseteq\iota(\tau)$,

	\item $\length{\sigma} =
		\length{\iota(\sigma)}$, and

	\item $\varphi(\sigma) =
		\varphi'(\iota(\sigma))$.
\end{enumerate}
The trees $\tree$ and $\tree'$ are
\emph{isomorphic}, denoted $\tree \isomorphic
\tree'$, if there exists an isomorphism between
them.

\begin{lem}
	\th\label{bisimilar_same_rank}
	If $B \subseteq \tree\times\tree'$ is a
	bisimulation and $\sigma \mathrel{B} \tau$
	holds, then $\rank_\tree(\sigma) =
	\rank_{\tree'}(\tau)$.
\end{lem}
\begin{proof}
	If $\rank_\tree(\sigma) = \infty$, i.e., if
	$\sigma$ is on an infinite branch of $\tree$,
	then it is easy to see that $\tau$ is on an
	infinite branch of $\tree'$ and therefore
	$\rank_{\tree'}(\tau) = \infty$ as well.
	By the same argument, we have that if
	$\rank_{\tree'}(\tau) = \infty$ then
	$\rank_\tree(\sigma) = \infty$.
	If $\sigma \in \WFpart(\tree)$, then we proceed
	by induction on $\rank_{\tree}(\sigma)$.
	For the base case, note that
	$\rank_{\tree}(\sigma) = 0$ iff $\sigma$ is a
	leaf of $\tree$, and in this case $\sigma
	\mathrel{B} \tau$
	implies that $\tau$ is also a leaf of $\tree'$
	and therefore also has rank $0$.
	Now suppose the result holds for every node of
	rank $< \rank_{\tree}(\sigma)$.
	For each $\beta < \rank_{\tree}(\sigma)$ there
	exists some descendant $\sigma'$ of $\sigma$ in
	$\tree$ such that
	$\rank_{\tree}(\sigma') = \beta$.
	Since $B$ is a bisimulation, there exists a
	descendant $\tau'$ of $\tau$ in $\tree'$ such
	that $\sigma' \mathrel{B} \tau'$.
	By induction hypothesis we have
	$\rank_{\tree'}(\tau') = \beta$, and since
	$\beta < \rank_{\tree}(\sigma)$ was arbitrary
	we have
	$\rank_{\tree'}(\tau) \geq
	\rank_{\tree}(\sigma)$.
	Analogously we can prove $\rank_{\tree'}(\tau)
	\leq \rank_{\tree}(\sigma)$, so the result
	follows.
\end{proof}

An \emph{abstract tree} is an equivalence class
of labeled trees under the relation of
bisimilarity.
Let $\AT$ be the space of abstract trees
represented by the total function $\delta_{\AT}$
given by $\delta_{\AT}(p) =
\delta_{\LT}(p)/{\bisim}$.
We typically denote abstract trees by
$\mathcal{A}$, with or without sub- or
superscripts.
As usual with quotient constructions, any
property of labeled trees can be extended to
abstract trees by stipulating that an abstract
tree has the property in question if one of its
representatives does.
Note that for some properties this extension
behaves better than for some others.
For example, the property of \emph{having rank
$\alpha$} behaves well, since by
\th\ref{bisimilar_same_rank} any two bisimilar
labeled
trees have the same rank.
On the other hand, the property of \emph{being
finitely branching} does not behave as well,
since every finitely branching labeled tree
is bisimilar to an infinitely branching one.

\begin{figure}[ht]
	\centering
	\includegraphics[scale=1]{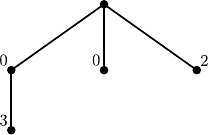}
	\caption{Depiction of the informal tree
	$\extset{(0,\extset{(3,\varnothing)}),(0,\varnothing),(2,\varnothing)}$.}
	\label{informal_tree}
\end{figure}

Note that, according to our definition, formally
speaking an abstract tree is not itself a tree
but only a certain type of set of
labeled trees.
However, for the sake of intuition it can be
helpful to think of an abstract tree as an
unordered tree without any concrete underlying
set of vertices, as follows.
We call an \emph{informal tree} a (possibly
empty) countable set $I$ of objects of the form
$(n,J)$, where $n$ is a natural number and
$J$ is again an informal tree.
The intuition is that such a tree $I$ is the tree
for which each such object $(n,J)$ represents a
child of the root of $I$ with label
$n$ and whose subtree is exactly $J$.
See Figure~\ref{informal_tree} for the depiction
of a simple informal tree.

To see how informal trees correspond to abstract
trees, let $\delta_\IT$ be the partial function
defined by corecursion with
\[
\delta_\IT(p)
= \intset{(n,\delta_\IT(q))}{\exists k ((p)_k =
\seq{n+1} \concat q)}.
\]
Then we say an informal tree $I$
\emph{corresponds} to an abstract tree
$\mathcal{A}$ if there exists $p \in
\dom{\delta_\IT}$ with
$\delta_\IT(p) = I$ and $\delta_\AT(p) =
\mathcal{A}$.

\begin{prop}
	In ZFC, the domain of $\delta_\IT$ is the set
	of $p \in \baire$ for which $\delta_\AT(p)$ is
	wellfounded.
	Therefore, in ZFC no informal tree corresponds
	to an illfounded abstract tree.
\end{prop}

This is, of course, because if $p$ is such that
$\delta_\AT(p)$ is illfounded, then in order for
$p \in \dom{\delta_\IT}$ to hold there
would have to exist an infinite $\in$-descending
chain of sets starting at $\delta_\IT(p)$,
contradicting the axiom of foundation.

However---as is often the case with definitions
by corecursion \cite{moss_danner_corecursion}---,
this definition and the correspondence
would also work for illfounded trees if one were
to work in a system of non-wellfounded set theory
such as ZFC\textsuperscript{$-$} +
AFA, where AFA is the axiom of anti-foundation
first formulated by Forti and
Honsell~\cite{forti_honsell_freeconstruction} and
later
popularized by
Aczel~\cite{aczel_nonwellfounded}---in the style
of Aczel~\cite[Chapter~6]{aczel_nonwellfounded},
in
ZFC\textsuperscript{$-$} + AFA the set of
informal trees can be defined as the greatest
fixed point of the class operator $\Phi$ defined
by letting $\Phi(X)$ be the class of all
countable sets of elements of the form $(n,T)$,
with $n \in \nat$ and $T$ a countable subset of
$X$.
Thus in ZFC\textsuperscript{$-$} + AFA the set of
informal trees is exactly
\[
\bigcup \intset{x \in V}{x \subseteq \Phi (x)}.
\]

\begin{prop}[ZFC\textsuperscript{$-$} + AFA]
	The correspondence between abstract and
	informal trees is a bijection.
\end{prop}

We will not pursue this line of investigation any
further; we thus now move back to our setting of
ZFC for the remainder of the paper.

\subsubsection*{Computable functions between
abstract trees}

We denote by $\open\nat$ the represented space of
subsets of $\nat$ given by \emph{enumeration},
i.e., so that $p$ is a name for $X
\subseteq \nat$ iff $X =
\intset{n\in\nat}{\exists k \in \nat(p(k) =
n+1)}$.
Note that any computable function of type
$\mpf{\open\nat}{\open\nat}$ has a computable
realizer which uses only \emph{positive
information}, by which we mean that it works by
only following rules of the form ``enumerate a
certain natural number into the output
set only after having seen some finite set of
natural numbers enumerated into the input set'',
i.e., via enumeration operators
(cf., e.g., \cite[Chapter
XIV]{odifreddi_book_volume_2}).
Let $\SubTrees:\function\LT{\open\nat}$ be
defined by letting $\SubTrees(\tree) = \intset{n
\in \nat}{\enumTrees(n)$ is a labeled
subtree of $\tree}$.
It is easy to see that $\SubTrees$ is computable.

\begin{lem}
	There exists a computable map $\ConstructTree:
	\mpf{\open\nat}{\LT}$ such that the composition
	$\ConstructTree \comp \SubTrees$ is total and
	$\tree' \in \ConstructTree \comp
	\SubTrees(\tree)$ implies $\tree' \isomorphic
	\tree$.
\end{lem}
\begin{proof}
	$\ConstructTree$ can be defined as follows.
	Suppose we are at stage $k$ of the
	construction, when some $n \in \nat$ is
	enumerated into the input.
	If some $m$ has been enumerated at some earlier
	stage such that $\enumTrees(m) \supset
	\enumTrees(n)$, then we proceed to the next
	stage.
	Otherwise let $X$ be the set of $m \in \nat$
	such that $\enumTrees(m)$ is a maximal subtree
	of $\enumTrees(n)$ among those $m$ which
	have been enumerated at earlier stages.
	By construction, for each $m \in X$ we have
	defined an associated $a(m) \in \nat$ and an
	isomorphism $\iota_m:
	\function{\enumTrees(m)}{\enumTrees(a(m))}$, in
	such a way that $\enumTrees(a(m))$ is
	guaranteed to be a subtree of the output tree
	we
	are constructing.
	Let $N \geq n$ be least such that there exists
	an isomorphism $\iota_n:
	\function{\enumTrees(n)}{\enumTrees(N)}$
	extending $\iota_m$
	for each $m \in X$ (in particular
	$\enumTrees(a(m)) \subset \enumTrees(N)$ for
	every $m \in X$) and such that no node of
	$\enumTrees(N)$ which is not in $\bigcup_{m\in
	X} \enumTrees(a(m))$ has been promised to be
	part of our current partial output.
	Then let $a(n):= N$ and guarantee that
	$\enumTrees(N)$ will be a subtree of our output
	tree.

	It is now straightforward to check that running
	the algorithm above on a name for
	$\SubTrees(\tree)$ we have $\tree =
	\bigcup_{n\in
	\dom{a}} \enumTrees(n)$ and that $\iota:=
	\bigcup_{n\in \dom{a}}\iota_n$ is an
	isomorphism between $\tree$ and $\tree':=
	\bigcup_{n\in \dom{a}} \enumTrees(a(n))$.
\end{proof}

\begin{lem}
	\th\label{modifying_realizers_AT}
	Let $G \realizes g:\mpf{\AT}{\AT}$.
	Suppose $F:\pf\baire\baire$ and
	$H:\pf\baire\baire$ are such that
	$\delta_{\LT}F(p) \bisim \delta_{\LT}(p)$ and
	$\delta_{\LT}H(q)
	\bisim \delta_{\LT}(q)$ for any $p \in \dom{F}$
	and $q \in \dom{H}$, and $\dom{G} \subseteq
	\dom{HGF}$.
	Then $HGF \realizes g$.
\end{lem}
\begin{proof}
	We have $\delta_{\AT} HGF(p) = \delta_{\AT}
	GF(p)$ and $\delta_{\AT}F(p) =
	\delta_{\AT}(p)$, therefore $\delta_{\AT}
	HGF(p) =
	\delta_{\AT}G(p)$.
\end{proof}

\begin{cor}
	\th\label{positive_information}
	Let $F, H$ be computable realizers of
	$\ConstructTree$ and $\SubTrees$, respectively.
	If $G$ is a computable realizer of some $g:
	\mpf{\AT}{\AT}$, then so is $FHGFH$.
\end{cor}
\begin{proof}
	Indeed, we have $\delta_{\LT}FH(p) \in
	\ConstructTree \comp \SubTrees
	\delta_{\LT}(p)$, so $\delta_{\LT}FH(p) \bisim
	\delta_{\LT}(p)$.
\end{proof}

Note that $HGF$ is a computable realizer of some
function $g': \mpf{\open\nat}{\open\nat}$; thus
we can assume that $FHGFH$ works by
only following rules of the form ``make a certain
finite labeled tree a subtree of the output only
after having seen some finite set of
finite labeled trees as subtrees of the input'',
which is to say, ``make a certain finite labeled
tree a subtree of the output only
after having seen a certain finite labeled tree
as a subtree of the input''.

\subsection{The pruning derivative}

We now define the main operation which will be
used in the game characterization of the class of
functions of any fixed Baire class.
First, let us recall the definition from
\cite{pauly_countable_ordinals} of the space
$\COrd$ of countable ordinals represented by the
function $\deltaCOrd$ defined recursively by
\begin{enumerate}
	\item $\deltaCOrd(0p) = 0$ \item
	$\deltaCOrd(1p) = \deltaCOrd(p) + 1$ \item
		$\deltaCOrd(2\tuple{p_n}_{n\in\nat}) =
		\sup_{n\in\nat}
		\deltaCOrd(p_n)$.
\end{enumerate}

\begin{defi}
	\th\label{pruning_derivative}
	We call \emph{pruning derivative} the operation
	$\PDName:\function{\UT}{\UT}$ which assigns to
	each unlabeled tree $T$ its subtree
	$\PD{T}:= \intset{\sigma \in T}{\sigma$ has
	descendants of arbitrary lengths in $T}$.
	We overload notation and also denote by
	$\PDName:\function{\LT}{\LT}$ the operation
	which assigns to each labeled tree $\tree =
	(T,\varphi)$ its subtree $\PD{\tree}$ whose
	domain is $\PD{T}$.
	As usual, these definitions can be iterated
	transfinitely in a natural way by letting
	\begin{enumerate}
		\item
			$\iPD{\argument}{0} = \id$
		\item
			$\iPD{\argument}{\alpha+1} =
			\PD{\iPD{\argument}{\alpha}}$
		\item $\iPD{\argument}{\lambda} =
			\bigcap_{\alpha<\lambda}
			\iPD{\argument}{\alpha}$, for limit
			$\lambda > 0$.
	\end{enumerate}
	Since trees are countable, for any given tree
	the iteration described above will stabilize at
	some countable stage for each tree.
	Thus, as functions between represented spaces,
	we can consider them as having types
	$\iPDName:\function{\UT \times \COrd}{\UT}$ and
	$\iPDName:\function{\LT \times \COrd}{\LT}$,
	respectively.
\end{defi}

\begin{lem}
	\th\label{derivative_rank}
	For $\sigma \in T$, we have $\sigma \in
	\iPD{T}{\alpha}$ iff $\rank_T(\sigma) \geq
	\omega\cdot \alpha$.
\end{lem}

\begin{lem}
	If $B \subseteq \tree\times\tree'$ is a
	bisimulation and $\sigma \mathrel{B} \tau$
	holds, then $\rank_\tree(\sigma) =
	\rank_{\tree'}(\tau)$.
\end{lem}
\begin{proof}
	By induction on $\rank_{\tree}(\sigma)$.
	For the base case, note that
	$\rank_{\tree}(\sigma) = 0$ iff $\sigma$ is a
	leaf of $\tree$, and in this case $\sigma
	\mathrel{B} \tau$
	implies that $\tau$ is also a leaf of $\tree'$
	and therefore also has rank $0$.
	Now suppose the result holds for every node of
	rank $< \rank_{\tree}(\sigma)$.
	For each $\beta < \rank_{\tree}(\sigma)$ there
	exists some descendant $\sigma'$ of $\sigma$ in
	$\tree$ such that
	$\rank_{\tree}(\sigma') = \beta$.
	Since $B$ is a bisimulation, there exists a
	descendant $\tau'$ of $\tau$ in $\tree'$ such
	that $\sigma' \mathrel{B} \tau'$.
	By induction hypothesis we have
	$\rank_{\tree'}(\tau') = \beta$, and since
	$\beta < \rank_{\tree}(\sigma)$ was arbitrary
	we have
	$\rank_{\tree'}(\tau) \geq
	\rank_{\tree}(\sigma)$.
	Analogously we can prove $\rank_{\tree'}(\tau)
	\leq \rank_{\tree}(\sigma)$, so the result
	follows.
\end{proof}

\begin{cor}
	\th\label{pruning_derivative_bisimilar}
	If $\tree \bisim \tree'$ then
	$\iPD{\tree}{\alpha} \bisim
	\iPD{\tree'}{\alpha}$ for any $\alpha <
	\omega_1$.
\end{cor}

We overload notation yet again and denote by
$\PDName:\function{\AT}{\AT}$ and
$\iPDName:\function{\AT\times\COrd}{\AT}$ the
operations
assigning to each abstract tree $\mathcal{A}$
with representative $\tree$ and each countable
ordinal $\alpha$ the subtrees
$\PD{\mathcal{A}}$ and
$\iPD{\mathcal{A}}{\alpha}$ with representatives
$\PD{\tree}$ and $\iPD{\tree}{\alpha}$,
respectively;
\th\ref{pruning_derivative_bisimilar} guarantees
that these are well-defined operations.
Whenever not specified otherwise by the context,
in what follows $\PDName$ and $\iPDName$ will
refer to the operations on abstract
trees.

\subsection{The Weihrauch degree of the pruning
derivative}

In order to analyze the Weihrauch degree of
$\iPDName$, we will first introduce and analyze
several operations on trees.
We introduce and analyze them as modularly as
possible, in the hope that this will increase the
clarity of the presentation and the
potential for applicability of the operations in
other situations.

Given $\sigma_0, \ldots, \sigma_{n-1} \in
\finbaire$ such that $\length{\sigma_i} =
\length{\sigma_j} = \ell$ for each $i,j < n$, let
$\finiteTreeProductNode{\sigma_0,\ldots,\sigma_{n-1}}
\in \nat^{\ell}$ be defined by
$\finiteTreeProductNode{\sigma_0,\ldots,\sigma_{n-1}}(m)
= \tuple{\sigma_0(m),\ldots,\sigma_{n-1}(m)}$ for
each $m < \ell$.
Note that
$\finiteTreeProductNode{\sigma_0,\ldots,\sigma_{n-1}}
\subseteq
\finiteTreeProductNode{\tau_0,\ldots,\tau_{n-1}}$
iff
$\sigma_i \subseteq \tau_i$ for every $i < n$.
Now, given trees $T_0,\ldots,T_{n-1}$, let their
\emph{product} be the tree
$\finiteTreeProduct{i<n}{T_i}:=
\intset{\finiteTreeProductNode{\sigma_0,\ldots,\sigma_{n-1}}}{\forall
i<n(\sigma_i \in T_i \text{ and }
\length{\sigma_{i}} = \length{\sigma_0})}$.
If $n = 2$ then we use the smaller infix notation
$\binaryTreeProduct{T_0}{T_1}$ to denote the
product.

\begin{lem}
	\th\label{tree_product}
	The operation
	$\finiteTreeProductName:\function{\finseq\UT}{\UT}$
	is computable and
	\begin{enumerate}
		\item $\finiteTreeProduct{i<n}{T_i} =
			\varnothing$ iff $T_i = \varnothing$ for
			some $i<n$.

		\item $\PD{\finiteTreeProduct{i<n}{T_i}} =
			\finiteTreeProduct{i<n}{\PD{T_i}}$.

		\item $\bigcap_{\beta < \alpha}
			\finiteTreeProduct{i<n}{T^\beta_i} =
			\finiteTreeProduct{i<n}{\bigcap_{\beta <
			\alpha} T^\beta_i}$
			for any ordinal $\alpha$.
	\end{enumerate}
\end{lem}

In particular,
$\iPD{\finiteTreeProduct{i<n}{T_i}}{\alpha} =
\finiteTreeProduct{i<n}{\iPD{T_i}{\alpha}}$ for
any ordinal $\alpha$.

We extend the binary product
$\binaryTreeProductName$ to type $\function{\LT
\times \UT}{\LT}$ by letting
$\binaryTreeProduct{(T,\varphi)}{S} =
(\binaryTreeProduct{T}{S},\varphi')$, where
$\varphi'(\tuple{\sigma,\tau}) =
\varphi(\sigma)$.

\begin{lem}
	\th\label{product_concrete_unlabeled_trees}
	If $S$ is pruned and nonempty, then
	$(T,\varphi) \bisim
	\binaryTreeProduct{(T,\varphi)}{S}$.
\end{lem}
\begin{proof}
	Let $B \subseteq T \times
	(\binaryTreeProduct{T}{S})$ be given by $\sigma
	\mathrel{B} \tau$ iff $\tau =
	\tuple{\sigma,\xi}$ for some
	$\xi \in S$.
	It is easy to see that $B$ satisfies conditions
	(\ref{lab}) and (\ref{par}).
	Suppose $\sigma \mathrel{B} \tau$, and let $\xi
	\in S$ be such that $\tau =
	\tuple{\sigma,\xi}$.
	For (\ref{for}), let $\sigma'$ be a child of
	$\sigma$ in $(T,\varphi)$.
	Since $S$ is pruned, $\xi$ has a child $\xi'$
	in $S$, and therefore $\tau':=
	\tuple{\sigma',\xi'}$ is a child of $\tau$ in
	$\binaryTreeProduct{(T,\phi)}{S}$.
	Now $\sigma' \mathrel{B} \tau'$ follows.
	For (\ref{bac}), let $\tau'$ be a child of
	$\tau$ in $\binaryTreeProduct{(T,\phi)}{S}$.
	Thus $\sigma'$ is a child of $\sigma$ in
	$(T,\varphi)$, from which $\sigma' \mathrel{B}
	\tau'$ follows.
\end{proof}

For our next operation on trees, let us first
define some auxiliary notation
$\infiniteTreeProductNode{\sigma_0,\ldots,\sigma_{n-1}}$,
for $\sigma_0,\ldots,\sigma_{n-1} \in \finbaire$.

\begin{defi}
	We define $\infiniteTreeProductNode{\,}:=
	\emptyseq$.
	Then, given $\sigma_0,\ldots,\sigma_{n-1} \in
	\finbaire$ such that $n = \length{\sigma_0} >
	0$ and $\length{\sigma_i} = n - i$ for
	each $i < n$, let
	$\infiniteTreeProductNode{\sigma_0,\ldots,\sigma_{n-1}}$
	be defined by
	\[
	\infiniteTreeProductNode{\sigma_0,\ldots,\sigma_{n-1}}(m)
	=
	\tuple{\sigma_0(m),\sigma_1(m-1),\sigma_2(m-2),\ldots,\sigma_m(0)}
	\]
	for each $m < n$.
	Note that
	$\infiniteTreeProductNode{\sigma_0,\ldots,\sigma_{n-1}}
	\subseteq
	\infiniteTreeProductNode{\tau_0,\ldots,\tau_{m-1}}$
	iff $n
	\leq m$ and $\sigma_i \subseteq \tau_i$ for
	each $i < n$.
	Now, given trees $\seq{T_n}_{n\in\nat}$, let
	their \emph{countable product} be the tree
	\[
	\infiniteTreeProduct{n\in\nat}{T_n}:=
	\intset{\infiniteTreeProductNode{
	\sigma_0,\ldots,\sigma_{k-1}}}{\forall m < k
	(\length{\sigma_m} = k-m \text{ and } \sigma_m
	\in T_m)}.
	\]
\end{defi}

Note that $\emptyseq \in
\infiniteTreeProduct{n\in\nat}{T_n}$ always
holds.
In particular, it is \emph{not} always the case
that
$\iPD{\infiniteTreeProduct{n\in\nat}{T_n}}{\alpha}
=
\infiniteTreeProduct{n\in\nat}{\iPD{T_n}{\alpha}}$
holds for all $\alpha$, contrary to the situation
for finite products of trees.

\begin{lem}
	\th\label{infinite_product}
	The operation
	$\infiniteTreeProductName:\function{\infseq\UT}{\UT}$
	is computable and
	\begin{enumerate}
		\item
			\label{infinite_product_rank}
			For $m > 0$ we have $\rank_S
			(\infiniteTreeProductNode{\sigma_0,\ldots,\sigma_{m-1}})
			\leq \min_{i<m} \rank_{T_i} (\sigma_i)$,
			with equality in case $\rank(T_j) \geq
			\min_{i<m} \rank_{T_i} (\sigma_i)$ holds
			for each $j \geq m$.
			As a consequence, we have $\rank(S) \leq
			\min_{n\in\nat} (\rank(T_n)+n)$.

		\item
			\label{infinite_product_rank_commutes_with_derivative}
			For every $\alpha$ we have
			$\iPD{\infiniteTreeProduct{n\in\nat}{T_n}}{\alpha}
			\subseteq
			\infiniteTreeProduct{n\in\nat}{\iPD{T_n}{\alpha}}$,
			with equality in case $\alpha = 0$ or
			$\iPD{T_n}{\alpha} \neq \varnothing$ for
			all $n \in \nat$.

		\item
			\label{infinite_product_nonempty_pruned}
			If all $T_n$ are pruned and nonempty then
			so is
			$\infiniteTreeProduct{n\in\nat}{T_n}$.
	\end{enumerate}
\end{lem}
\begin{proof}
	The computability of $\infiniteTreeProductName$
	is straightforward.

	(\ref{infinite_product_rank})
	By induction on $\rank_{T_i} (\sigma_i)$, we
	show that $\rank_S
	(\infiniteTreeProductNode{\sigma_0,\ldots,\sigma_{m-1}})
	\leq
	\rank_{T_i} (\sigma_i)$.
	If $\rank_{T_i} (\sigma_i) = 0$, this is easy
	to see.
	For $\rank_{T_i} (\sigma_i) > 0$ we have that
	every descendant of $\sigma_i$ in $T_i$ has
	rank less than $\rank_{T_i} (\sigma_i)$, so
	by inductive hypothesis every descendant of
	$\infiniteTreeProductNode{\sigma_0,\ldots,\sigma_{m-1}}$
	in $S$ has rank less than
	$\rank_{T_i} (\sigma_i)$, and therefore
	$\rank_S
	(\infiniteTreeProductNode{\sigma_0,\ldots,\sigma_{m-1}})
	\leq \rank_{T_i}
	(\sigma_i)$.
	Conversely, by induction on $\alpha$ we show
	that if $\rank_{T_i} (\sigma_i), \rank(T_j)
	\geq \alpha$ holds for all $i < m$ and $j
	\geq m$, then $\rank_S
	(\infiniteTreeProductNode{\sigma_0,\ldots,\sigma_{m-1}})
	\geq \alpha$ as well.
	The case $\alpha = 0$ is clear.
	Now suppose $\alpha > 0$.
	Given $\beta < \alpha$, for each $i < m$ let
	$\sigma'_i$ be an immediate child of $\sigma_i$
	in $T_i$ of rank at least $\beta$, and
	let $\sigma'_m \in T_m$ have length $1$ and
	rank at least $\beta$.
	Then
	$\infiniteTreeProductNode{\sigma'_0,\ldots,\sigma'_m}$
	is an immediate child of
	$\infiniteTreeProductNode{\sigma_0,\ldots,\sigma_{m-1}}$
	in $S$.
	Since $\rank(T_j) \geq \beta$ for each $j \geq
	m+1$, by induction hypothesis we get $\rank_S
	(\infiniteTreeProductNode{\sigma'_0,\ldots,\sigma'_m})
	\geq \beta$.
	Therefore $\rank_S
	(\infiniteTreeProductNode{\sigma_0,\ldots,\sigma_{m-1}})
	> \beta$, and since $\beta < \alpha$ was
	arbitrary we get
	$\rank_S
	(\infiniteTreeProductNode{\sigma_0,\ldots,\sigma_{m-1}})
	\geq \alpha$, as desired.
	Finally $\tau \in S$ has length $n+1$, then
	$\tau =
	\infiniteTreeProductNode{\sigma_0,\ldots,\sigma_{n}}$
	where $\sigma_i \in T_i$ for
	every $i \leq n$.
	In particular, $\rank_S (\tau) \leq \rank_{T_n}
	(\sigma_n) < \rank(T_n)$, so $\rank(S) =
	\rank_S(\emptyseq) \leq \rank(T_n) + n$.

	(\ref{infinite_product_rank_commutes_with_derivative})
 Follows by combining
 (\ref{infinite_product_rank}) with
 \th\ref{derivative_rank}.

 (\ref{infinite_product_nonempty_pruned})
 Follows from (\ref{infinite_product_rank}) since
 a tree is pruned and nonempty iff all its nodes
 have rank~$\infty$.
\end{proof}

\begin{defi}
	Given trees $T_0,\ldots,T_{n-1}$, let their
	\emph{mix} be the tree $\finiteMix{i<n}{T_i}$
	which satisfies $\finiteMix{i<n}{T_i} =
	\varnothing$ iff $T_i = \varnothing$ for some
	$i < n$, and otherwise $\emptyseq \in
	\finiteMix{i<n}{T_i}$ and
	$\concatenationSubtree{\finiteMix{i<n}{T_i}}{\seq{\tuple{m,k}}}
	= \concatenationSubtree{T_m}{\seq{k}}$ for each
	$m < n$ and $k \in
	\nat$.
	Intuitively, the mix of $T_0,\ldots,T_{n-1}$ is
	the tree obtained by merging the roots of those
	trees into a single root.
	If $n = 2$ then we use the smaller infix
	notation $\binaryMix{T_0}{T_1}$ to denote the
	mix.
\end{defi}

\begin{lem}
	\th\label{mix}
	The operation
	$\finiteMixName:\function{\finseq\UT}{\UT}$ is
	computable and
	\begin{enumerate}
		\item $\PD{\finiteMix{i<n}{T_i}} =
			\finiteMix{i<n}{\PD{T_i}}$.

		\item $\bigcap_{\beta < \alpha}
			\finiteMix{i<n}{T^\beta_i} =
			\finiteMix{i<n}{\bigcap_{\beta < \alpha}
			T^\beta_i}$ for any ordinal
			$\alpha$.
	\end{enumerate}

	In particular,
	$\iPD{\finiteMix{i<n}{T_i}}{\alpha} =
	\finiteMix{i<n}{\iPD{T_i}{\alpha}}$ for any
ordinal $\alpha$ \end{lem}

\begin{defi}
	Given trees $\seq{T_n}_{n\in\nat}$, let their
	\emph{countable mix} be the tree
	$\infiniteMix{n\in\nat}{T_n}$ such that
	$\emptyseq \in
	\infiniteMix{n\in\nat}{T_n}$ and
	$\concatenationSubtree{\infiniteMix{n\in\nat}{T_n}}{\seq{\tuple{m,k}}}
	=
	\concatenationSubtree{T_m}{\seq{k}}$ for each
	$m,k \in \nat$.
	As before with countable products,
	$\infiniteMixName$ will not commute with
	$\iPD{\argument}{\alpha}$ for all $\alpha$ in
	general.
\end{defi}

\begin{lem}
	\th\label{infinite_mix}
	The operation
	$\infiniteMixName:\function{\infseq\UT}{\UT}$
	is computable and satisfies
	$\iPD{\infiniteMix{n\in\nat}{T_n}}{\alpha}
	\subseteq
	\infiniteMix{n\in\nat}{\iPD{T_n}{\alpha}}$,
	with equality in case $\alpha = 0$ or
	$\iPD{T_n}{\alpha} \neq \varnothing$ for
	some $n \in \nat$.
\end{lem}

To proceed, we need the notion of a Borel truth
value.
This represented space was introduced in
\cite{gregoriades_kispeter_pauly_comparison}
(built on ideas from \cite{moschovakis_book}),
and
further investigated in
\cite{pauly_countable_ordinals}.
Our definition differs slightly from the one
given in the literature, but is easily seen to be
equivalent.

\begin{defi}
	A \emph{Borel truth value} is a pair $b =
	(T,\mu)$ such that $T$ is a wellfounded tree
	and $\mu$ is a function, called a
	\emph{tagging} function, assigning to each node
	of $T$ one of the \emph{tags}
	$\bot,\top,\forall,\exists$, in such a way that
	each
	leaf is tagged $\top$ or $\bot$, and each
	non-leaf node is tagged $\forall$ or $\exists$
	(in alternating fashion, i.e., so that if a
	node tagged $\forall$ has a parent, then the
	parent is tagged $\exists$ and vice versa).
	A name for a Borel truth value $(T,\mu)$ is an
	element $p \in \infseq{5}$ which is a
	$\delta_\UT$ -code for $T$ and such that if
	$\sigma \in T$, i.e., if $\sigma$ is a path
	through $p$, then
	\[
	\mu(\sigma) =
	\begin{cases}
		\bot, &\text{if } \sigma = \emptyseq \text{
		and } p(0)=1 \text{, or } \length\sigma = 1
		\text{ and } p(\sigma(0))=1 \\
		\top, &\text{if } \sigma = \emptyseq \text{
		and } p(0)=2 \text{, or } \length\sigma = 1
		\text{ and } p(\sigma(0))=2 \\
		\forall, &\text{if } \sigma = \emptyseq
		\text{ and } p(0)=3 \text{, or }
		\length\sigma = 1 \text{ and } p(\sigma(0))=3
		\\
		\exists, &\text{if } \sigma = \emptyseq
		\text{ and } p(0)=4 \text{, or }
		\length\sigma = 1 \text{ and } p(\sigma(0))=4
		\\
		\mu'(\shift{\sigma}), &\text{if }
		\length\sigma > 1,
	\end{cases}
	\]
	where $(T',\mu')$ is the Borel truth value
	named by $(\shift{p})_{\sigma(0)}$.
	In other words, intuitively in a name for a
	Borel truth value, zeroes indicate absence of
	the corresponding node, and nonzero values
	indicate both presence of the corresponding
	node and its tag.
	The \emph{value} $\val(b) \in
	\extset{\top,\bot}$ of a Borel truth value $b$
	is defined by recursion on the rank of $b$ in a
	straightforward way.
	The space of Borel truth values is denoted by
	$\mathbb{S}(\mathcal{B})$.
	The \emph{$\Sigma^0_\alpha$ -truth values}
	(denoted $\mathbb{S}({\Sigma^0_\alpha})$) are
	those with rank $\leq \alpha$ and root tagged
	$\exists$, and the \emph{$\Pi^0_\alpha$ -truth
	values} (denoted $\mathbb{S}({\Pi^0_\alpha})$)
	are those with rank $\leq \alpha$ and
	root tagged $\forall$.
\end{defi}

Given an ordinal $\alpha$, with $\alpha = \lambda
+ n$ for some limit ordinal $\lambda$ and $n \in
\nat$, let $\hat{\alpha} = \lambda +
2n$ and $\check{\alpha} = \lambda + \lceil
\frac{n}{2} \rceil$.

\begin{prop}
	\th\label{prop:ispresent}
	The map $\mathrm{isPresent}: \mpf{\LT \times
	\COrd \times \nat}{\coprod_{\alpha \in \COrd}
	\mathbb{S}({\Pi^0_{\alpha}}})$, mapping
	$(\tree,\alpha,\ell)$ such that
	$\enumTrees(\ell)$ is linear to $(\max
	\extset{1,\hat{\alpha}}, b)$ where $\val(b) =
	\top$ iff
	$\enumTrees(\ell) \subseteq
	\iPD{\tree}{\alpha}$, is computable.
\end{prop}
\begin{proof}
	It is straightforward to see that $\alpha
	\mapsto \max \extset{1,\hat\alpha}: \COrd \to
	\COrd$ is computable.
	Computability of the second component is shown
	by induction over the $\deltaCOrd$ -name $q$ of
	$\alpha$ provided.
	If $q = 0q'$, then we check whether
	$\enumTrees(\ell) \subseteq \tree$ and return
	either the tree of rank $1$ with root tagged
	$\forall$ and children tagged $\top$ (if yes),
	or with root tagged $\forall$ and children
	tagged $\bot$ (if no).
	If $q = 1q'$, then $\alpha = \beta + 1$ and
	$q'$ is a name for $\beta$.
	Let $h$ be the height of $\enumTrees(\ell)$.
	We start searching for confirmation that $\beta
	> 0$.
	Until we find it, we add children with tag
	$\exists$ to the root tagged $\forall$, and
	then for each $\ell' \in \nat$ such that
	$\enumTrees(\ell')$ is a linear tree of height
	$h'>h$ extending $\enumTrees(\ell)$, we add a
	grandchild tagged $\top$ or $\bot$ to the
	$(h'-h)$ \textsuperscript{th} child, depending
	on whether or not $\enumTrees(\ell') \subseteq
	\tree$.
	If we do receive confirmation that $\beta > 0$,
	we add a grandchild tagged $\top$ to each
	$\exists$ -child produced so far, and then
	ignore these children.
	Then, for each $\ell' \in \nat$ such that
	$\enumTrees(\ell')$ is a linear tree of height
	$h'>h$ extending $\enumTrees(\ell)$, we
	compute $b_{\ell'} \in \mathbb{S}({\Pi^0_{\max
	\extset{1,\hat\beta}}})$ denoting whether or
	not $\enumTrees(\ell') \subseteq
	\iPD{\tree}{\beta}$.
	Then we add each $b_{\ell'}$ as a grandchild of
	the root of $b$ via the $(h'-h)$
	\textsuperscript{th} new child tagged
	$\exists$.
	If $q = 2\tuple{q_i}_{i\in\nat}$, then $\alpha
	= \sup_{i \in \nat} \alpha_i$ and each $q_i$ is
	a name for $\alpha_i$.
	For each $i \in \nat$, we compute whether
	$\enumTrees(\ell) \subseteq
	\iPD{\tree}{\alpha_i}$ as $b_i \in
	\mathbb{S}({\Pi^0_{\max
	\extset{1,\hat{\alpha_i}}}})$.
	By induction, each $b_i$ has root tagged
	$\forall$, and we now obtain the answer $b$ as
	the mix of the $b_i$.

	\begin{clm}
		$\val(b) = \top$ iff $\enumTrees(\ell)
		\subseteq \iPD{\tree}{\alpha}$.
	\end{clm}
	By induction on the name $q$ of $\alpha$.
	If $q = 0q'$ then it is immediate to see that
	the claim holds.
 Suppose the claim holds for $q'$ and let $q=
 1q'$.
 Let $\beta = \deltaCOrd(q')$ and suppose
 $\enumTrees(\ell)$ has height $h$.
 Then $\enumTrees(\ell) \subseteq
 \iPD{\tree}{\alpha}$ iff for every $n \in \nat$
 there exists some $\ell' \in \nat$ such that
 $\enumTrees(\ell')$ is a linear tree of height
 $h+n$ extending $\enumTrees(\ell)$, with
 $\enumTrees(\ell') \subseteq
 \iPD{\tree}{\beta}$.
 By induction, the result of iterating the
 algorithm for $(\tree,\beta,\ell')$ gives the
 correct output.
 Therefore $\val(b) = \top$ iff $\enumTrees(\ell)
 \subseteq \iPD{\tree}{\alpha}$.
 Finally, suppose the claim holds for each $q_n$
 and let $q = 2\tuple{q_n}_{n\in\nat}$.
 Let $\alpha_n = \deltaCOrd(q_n)$.
 Again, by induction the result of iterating the
 algorithm for $(\tree,\alpha_n,\ell)$ gives the
 correct output.
 Therefore $\val(b) = \top$ iff all children of
 the roots of all $b_n$ have value $\top$ iff
 $\enumTrees(\ell) \subseteq
 \iPD{\tree}{\alpha_n}$ for all $n \in \nat$ iff
 $\enumTrees(\ell) \subseteq
 \iPD{\tree}{\alpha}$, as desired.

 \begin{clm}
	 The Borel truth value $b$ has rank $\leq \max
	 \extset{1,\hat\alpha}$.
 \end{clm}

 We again proceed by induction on the name $q$ of
 $\alpha$.
 If $q = 0q'$ then by construction $b$ has rank
 $1$.
 If $q = 1q'$ and $\beta:= \deltaCOrd(q') = 0$,
 then again by construction $b$ has rank $\leq 2
 = \hat{1}$.
 If $q = 1q'$ and $\beta:= \deltaCOrd(q') > 0$,
 then by induction each $b_\tau$ as defined in
 the algorithm has rank $\leq \hat\beta$,
 and therefore $b$ has rank $\leq \hat\beta + 2 =
 \hat\alpha$.
 Finally, if $q = 2\tuple{q_i}_{i\in\nat}$, for
 each $i$ let $\alpha_i = \deltaCOrd(q_i)$.
 By induction, each $b_i$ as defined in the
 algorithm has rank $\leq \max
 \extset{1,\hat{\alpha_i}}$, and by construction
 $b$ has rank
 $\leq \max \extset{1,\sup_{i\in\nat}
 \hat{\alpha_i}} = \max \extset{1,\hat\alpha}$.
\end{proof}

\begin{prop}
	\th\label{prop:witness}
	The map $\mathrm{Witness}:
	\mpf{\mathbb{S}(\mathcal{B})}{\UT}$, mapping
	$b$ of rank $\alpha>0$ to some $T$ such that if
	$\val(b) =
	\top$ then $\iPD{T}{\check{\alpha}}$ is a
	nonempty pruned tree, and if $\val(b) = \bot$
	then $\iPD{T}{\check{\alpha}} = \varnothing$,
	is computable.
\end{prop}
\begin{proof}
	If $b$ is composed of a single node, then we
	output $\finbaire$ or $\varnothing$ according
	to whether $\val(b) = \top$ or $\val(b) =
	\bot$.
	Otherwise, we iteratively compute trees
	$(T_n)_{n \in \nat}$ for all the subtrees
	rooted at the children of the root of $b$, and
	output $\infiniteTreeProduct{n \in \nat}{T_n}$
	if the root of $b$ is tagged $\forall$, or
	output $\infiniteMix{n \in \nat}{T_n}$ if
	the root of $b$ is tagged $\exists$.

	\begin{clm}
		Suppose $\beta > 0$ is such that
		$\iPD{T_n}{\beta}$ is pruned for each $n \in
		\nat$.
		Then $\iPD{T}{\beta}$ is pruned, if the root
		of $b$ is tagged $\forall$, and
		$\iPD{T}{\beta+1}$ is pruned, if the root of
		$b$ is
		tagged $\exists$.
		Furthermore, if $\iPD{T}{\delta}$ is pruned,
		then it is nonempty in case $\val(b) = \top$
		and empty in case $\val(b) = \bot$.
	\end{clm}

	If $b$ is composed of a single node then the
	claim follows easily.
	Otherwise, suppose the root of $b$ is tagged
	$\forall$, so that $T =
	\infiniteTreeProduct{n\in\nat}{T_n}$.
	If $\val(b) = \top$ then each
	$\iPD{T_n}{\beta}$ is pruned and nonempty, and
	therefore the same holds for $\iPD{T}{\beta}$.
	Conversely, if $\val(b) = \bot$ then
	$\iPD{T_{n_0}}{\beta} = \varnothing$ for some
	$n_0 \in \nat$.
	Thus there is some $\gamma < \beta$ and $H \in
	\nat$ such that $\iPD{T_{n_0}}{\gamma}$ has
	height $\leq H < \omega$.
	Then $\iPD{T}{\gamma}$ has height $\leq H' <
	\omega$ for some $H'$ depending on $H$ and
	$n_0$, and therefore $\iPD{T}{\beta} =
	\varnothing$.
	Now suppose the root of $b$ is tagged
	$\exists$, so that $T =
	\infiniteMix{n\in\nat}{T_n}$.
	If $\val(b) = \top$ then some
	$\iPD{T_n}{\beta}$ is pruned and nonempty.
	Therefore the same holds for $\iPD{T}{\beta}$.
	Otherwise, if $\val(b) = \bot$ then each
	$\iPD{T_n}{\beta}$ is empty.
	Therefore $\iPD{T}{\beta} \subseteq
	\extset{\emptyseq}$, and thus
	$\iPD{T}{\beta+1}$ is empty.

	\begin{clm}
		Let $b'$ be a Borel truth value and let $T'$
		be the result of applying the algorithm above
		to $b'$.
		For $\beta = \rank(b')$, we have that if the
		root of $b'$ has tag $\forall$ then
		$\iPD{T'}{\check\beta}$ is pruned, and if the
		root
		of $b'$ has tag $\exists$ then
		$\iPD{T'}{\check\beta+1}$ is pruned.
	\end{clm}
 By induction on $\beta$.
 If $\beta = 0$, i.e., if $b'$ is a single node,
 then by construction $T'$ is pruned.
 Now suppose $\beta > 0$, and let the $n$
 \textsuperscript{th} child $\sigma_n$ of the
 root of $b'$ have rank $\beta_n < \beta$.
 Suppose the root of $b'$ has tag $\exists$, so
 that each $\sigma_n$ is either a leaf or has tag
 $\forall$.
 By induction, the result $T_n$ of applying the
 algorithm to the subtree of $b'$ rooted at
 $\sigma_n$ is such that
 $\iPD{T_n}{\check\beta_n}$ is pruned.
 Since $\sup_{n\in\nat} \check\beta_n \leq
 \check\beta$, by the preceding claim it follows
 that $\iPD{T}{\check\beta+1}$ is pruned.
 Finally, suppose the root of $b'$ has tag
 $\forall$, so that each $\sigma_n$ is either a
 leaf or has tag $\exists$.
 By induction, the result $T_n$ of applying the
 algorithm to the subtree of $b'$ rooted at
 $\sigma_n$ is such that
 $\iPD{T_n}{\check\beta_n+1}$ is pruned.
 If $\sup_{n\in\nat}(\check\beta_n+1) \leq
 \check\beta$ for each $n$, then by the preceding
 claim we are done.
 Otherwise, say $\check\beta_n = \check\beta$ for
 some $n$.
 Then $\beta_n$ is odd and $\beta = \beta_n+1$.
 In particular $\beta_n = \gamma_n + 1$ for some
 $\gamma_n$, and therefore $\check\delta \leq
 \check\gamma_n < \check\beta_n$ for each
 $\delta < \beta_n$.
 Hence, since by induction the result $S$ of
 applying the algorithm to any subtree of
 $\beta'$ rooted at some child of $\sigma_n$ is
 such that $\iPD{S}{\check\delta}$ is pruned for
 some $\delta < \beta_n$, by the preceding claim
 it follows that
 $\iPD{T_n}{\check\gamma_n+1} =
 \iPD{T_n}{\check\beta_n}$ is pruned.
 Thus, again by the preceding claim,
 $\iPD{T}{\check\beta}$ is pruned.
\end{proof}

\begin{lem}
	For each $\alpha \in \COrd$ we have that
	$\iPD{\argument}{\alpha}$ is parallelizable.
\end{lem}
\begin{proof}
	Given abstract trees
	$\seq{\mathcal{A}_n}_{n\in\nat}$ with
	respective representatives
	$\seq{\tree_n}_{n\in\nat}$, let
	$\mathcal{A}$ be the abstract tree represented
	by the labeled tree $\tree$ in which the root
	has a child $\sigma_n$ labeled $n$ for
	each $n \in \nat$ such that $\tree_n \neq
	\varnothing$, and such that
	$\concatenationSubtree{\tree}{\sigma_n} =
	\tree_n$ in the
	positive case.
	It is now straightforward to see that each
	$\iPD{\mathcal{A}_n}{\alpha}$ can be
	reconstructed from $\iPD{\mathcal{A}}{\alpha}$.
\end{proof}

Let $\twoSpace$ be the represented space composed
of two elements, $\top$ and $\bot$, the first
represented by $1^\nat$ and the latter
by $0^\nat$.

\begin{cor}
	For each $\alpha > 0$ we have that
	$\iPD{\argument}{\alpha}$ is
	Weihrauch-equivalent to the parallelization of
	$\id_\alpha:
	\mathbb{S}({\Pi^0_{\hat{\alpha}}}) \to
	\twoSpace$.
	Furthermore, the reductions in both directions
	can be taken to be uniform in $\alpha$.
\end{cor}
\begin{proof}
	To reduce $\iPD{\argument}{\alpha}$ to the
	parallelization of $\id_\alpha$, note that we
	can use $\mathrm{isPresent}$ from
	\th\ref{prop:ispresent} to compute for each
	linear labeled tree whether or not to include
	it in $\iPD{\tree}{\alpha}$ as a
	$\mathbb{S}({\Pi^0_{\hat{\alpha}}})$-truth
	value.
	We then use the parallelization of $\id_\alpha:
	\mathbb{S}({\Pi^0_{\hat{\alpha}}}) \to
	\twoSpace$ to convert all of these into
	booleans, and can thus construct
	$\iPD{\tree}{\alpha}$.
	For the converse, we use
	$\mathrm{Witness}(\hat{\alpha},b)$ from
	\th\ref{prop:witness} to obtain some $\tree$
	such that
	$\iPD{\tree}{\alpha} = \varnothing$ if $\val(b)
	= \bot$ and $\iPD{\tree}{\alpha} \neq
	\varnothing$ if $\val(b) = \top$.
	As $\extset{\varnothing}$ is a decidable subset
	of $\LT$, we can recover $b \in \twoSpace$
	after obtaining $\iPD{\tree}{\alpha}$ from
	the oracle.
\end{proof}

\begin{thm}[Folklore]
	If $A \subseteq \baire$ is Wadge-complete for
	$\borelPi{\alpha}$ (respectively effectively
	Wadge-complete for $\Pi^0_\alpha$), then
	$\para{\chi_A}$ is continuously
	Weihrauch-complete for Baire class $\alpha$
	(respectively Weihrauch-complete for effective
	Baire class
	$\alpha$), where
	$\chi_A:\function{\baire}{\twoSpace}$ is given
	by $\chi_A(x) = \top$ iff $x \in A$.
\end{thm}
\begin{proof}
	That $\para{\characteristicFunction{A}}$ is
	Baire class $\alpha$ follows from noticing that
	\[
	\preimage{\para{\characteristicFunction{A}}}[\sigma]
	=
	\bigcap_{\substack{n < \length{\sigma} \\
	\sigma(n) = 0}} \intset{x \in \baire}{(x)_n \in
	A} \cap \bigcap_{\substack{n <
	\length{\sigma} \\ \sigma(n) = 1}} \intset{x
	\in \baire}{(x)_n \not\in A},
	\]
	which is the intersection of a
	$\borelPi{\alpha}$ set with a
	$\borelSigma{\alpha}$ set.

	Now let $F:\pf\baire\baire$ be a Baire class
	$\alpha$ realizer of
	$f:\mpf{\repSpace{X}}{\repSpace{Y}}$.
	Let $\intseq{\sigma_n}{n\in\nat}$ be some
	enumeration of $\finbaire$.
	Since $F$ is Baire class $\alpha$, there exists
	some countable collection
	$\intseq{X_{n,m}}{n,m\in\nat}$ of
	$\borelSigma{\alpha}$ sets
	such that $\preimage{F}[\sigma_n] =
	\bigcup_{m\in\nat} X_{n,m}$.
	Since $A$ is Wadge-complete for
	$\borelPi{\alpha}$, for each $n,m \in \nat$
	there exists a continuous
	$f_{n,m}:\function{\baire}{\baire}$ such that
	$X_{n,m} = \preimage{f_{n,m}}[\baire
	\smallsetminus A]$.
	Now, defining a continuous
	$K:\function{\baire}{\baire}$ by
	$(K(x))_{\tuple{n,m}} = f_{n,m}(x)$, we have
	$\sigma_n \subseteq F(x)$ iff
	$x \in X_{n,m}$ for some $m$ iff
	$\para{\characteristicFunction{A}}(K(x))(\tuple{n,m})
	= 1$ for some $m$.
	Finally, defining a continuous
	$H:\pf\baire\baire$ by $H(x) = \bigcup
	\intset{\sigma_n}{\exists m(x(\tuple{n,m}) =
	1}$ with its
	natural domain, we have
	$H\para{\characteristicFunction{A}}K \tightens
	F$.
	Therefore $F \leqsW
	\para{\characteristicFunction{A}}$, and $f
	\leqsW \para{\characteristicFunction{A}}$ as
	well.
\end{proof}

\begin{cor}
	For each $\alpha > 0$ the parallelization of
	the map $\id_\alpha:
	\mathbb{S}({\Pi^0_{\hat{\alpha}}}) \to
	\twoSpace$ is (continuously)
	Weihrauch-complete for (effective) Baire class
	$\alpha$.
\end{cor}
\begin{proof}
	It is enough to show that for each $\alpha > 0$
	the characteristic function of any
	$\borelPi{\alpha}$ set is continuously
	Weihrauch-reducible to $\id_\alpha$, and that
	$\id_\alpha$ is Weihrauch-reducible to the
	characteristic function of some
	$\borelPi{\alpha}$ set.
	Both of these claims can be easily proved by
	induction.
\end{proof}

\begin{cor}
	$\iPD{\argument}{\alpha}$ is Weihrauch-complete
	for Baire class $\hat\alpha$.
\end{cor}

\begin{cor}
	$\iPD{\argument}{\argument} \equivW
	\UC_\Baire$.
\end{cor}
\begin{proof}
	It was shown in \cite{pauly_countable_ordinals}
	that $\widehat{\id: \mathbb{S}(\mathcal{B}) \to
	\twoSpace} \equivW \UC_\Baire$.
\end{proof}

\subsection{The transparency of the pruning
derivative}

\begin{prop}
	\th\label{prop:isabsent}
	The map $\mathrm{isAbsent}: \mpf{\LT \times
	\COrd \times \nat}{\coprod_{\alpha \in \COrd}
	\mathbb{S}({\Pi^0_{\alpha}}})$, mapping
	$(\tree,\alpha,\ell)$ such that
	$\enumTrees(\ell)$ is linear to $(\max
	\extset{1,\ordinalup{\alpha}}, b)$ where
	$\val(b) = \top$ iff
	$\enumTrees(\ell) \not\subseteq
	\iPD{\tree}{\alpha}$, is computable.
\end{prop}
\begin{proof}
	We just run $\mathrm{isPresent}$ on
	$(\tree,\alpha,\ell)$ then dualize the output
	by exchanging tags $\forall$ with $\exists$ and
	$\top$ with $\bot$.
\end{proof}

\begin{cor}
	\th\label{witnessleq}
	The operation
	$\WitnessleqName:\mf{\LT\times\COrd\times\nat}{\UT}$,
	given by $S \in \Witnessleq{\tree,\alpha,m}$
	iff $\iPD{S}{\alpha}$
	is a pruned tree and $\iPD{S}{\alpha} \neq
	\varnothing$ iff $\enumTrees(m) \subseteq
	\iPD{\tree}{\alpha}$, is computable.
\end{cor}
\begin{proof}
	Let $\tree_0,\ldots,\tree_{n}$ be the linear
	subtrees of $\enumTrees(m)$, and for each $i
	\leq n$ let $\ell_{i} \in \nat$ be such that
	$\enumTrees(\ell_i) = \tree_i$.
	By \th\ref{prop:ispresent,prop:witness},
	letting $S_i \in \mathrm{Witness} \comp
	\mathrm{isPresent}(\tree,\alpha,\ell_i)$, we
	have that
	$\iPD{S_i}{\alpha}$ is pruned and
	$\iPD{S_i}{\alpha} \neq \varnothing$ iff
	$\tree_i \subseteq \iPD{\tree}{\alpha}$.
	Now letting $S = \finiteTreeProduct{i\leq
	n}{S_i}$ we have that $\iPD{S}{\alpha}$ is
	pruned and that $\iPD{S}{\alpha} \neq
	\varnothing$ iff $\enumTrees(m) \subseteq
	\iPD{\tree}{\alpha}$, as desired.
\end{proof}

\begin{cor}
	\th\label{empty}
	The operation
	$\WitnessAbsence:\mf{\LT\times\COrd\times\baire}{\UT}$,
	given by $S \in
	\WitnessAbsence(\tree,\alpha,x)$ iff in case
	$\alpha > 0$ then $\iPD{S}{\alpha}$ is a pruned
	tree and $\iPD{S}{\alpha} = \varnothing$ iff
	$\enumTrees(x(m)) \subseteq
	\iPD{\tree}{\alpha}$ for some $m \in \nat$, is
	computable.
\end{cor}
\begin{proof}
	For each $m$ and each linear subtree
	$\tree^m_0,\ldots,\tree^m_{n_m}$ of
	$\enumTrees(x(m))$, let $\ell^m_{i} \in \nat$
	be such that
	$\enumTrees(\ell^m_i) = \tree^m_i$.
	By \th\ref{prop:isabsent,prop:witness}, letting
	$S^m_i \in
	\mathrm{Witness}\comp\mathrm{isAbsent}(\tree,\alpha,\ell^m_i)$,
	we have
	that $\iPD{S^m_i}{\alpha}$ is a pruned tree and
	$\iPD{S^m_i}{\alpha} = \varnothing$ iff
	$\tree^m_i \subseteq \iPD{\tree}{\alpha}$.
	Now letting $S^m = \finiteTreeProduct{i\leq
	n_m}{S^m_i}$ we have that $\iPD{S^m}{\alpha}$
	is pruned and that $\iPD{S^m}{\alpha} =
	\varnothing$ iff $\enumTrees(x(m)) \subseteq
	\iPD{\tree}{\alpha}$.
	Now let $S = \infiniteTreeProduct{m \in
	\nat}{S^m}$.
	Suppose $\alpha > 0$, and first suppose that
	$\enumTrees(x(m)) \subseteq
	\iPD{\tree}{\alpha}$ for some $m \in \nat$.
	Then $\iPD{S^m}{\alpha} = \varnothing$, so for
	some $\beta < \alpha$ we have that
	$\iPD{S^m}{\beta}$ has some finite height $H$.
	Hence $\iPD{S}{\beta}$ also has some finite
	height $H'$ (which depends on $H$ and $m$), and
	therefore $\iPD{S}{\alpha} = \varnothing$,
	as desired.
	Now suppose $\enumTrees(x(m)) \not\subseteq
	\iPD{\tree}{\alpha}$ for all $m \in \nat$.
	Then each $\iPD{S^m}{\alpha}$ is pruned and
	nonempty, and therefore the same holds for
	$\iPD{S}{\alpha}$.
\end{proof}

\begin{prop}[Pauly {\cite[Theorem
	31]{pauly_countable_ordinals}}]
	\th\label{min_computable}
	The function
	$\min:\function{\COrd\times\COrd}{\COrd}$ is
	computable.
\end{prop}

\begin{prop}
	The function $\TreeWithRank:\mf{\COrd}{\UT}$,
	given by
	\[
	T \in \TreeWithRank(\alpha)
	\;\;\;\text{iff}\;\;\; T \text{ is a
	wellfounded tree and } \rank(T) = \alpha,
	\]
	is computable.
\end{prop}
\begin{proof}
	We will define a computable $F \realizes
	\TreeWithRank$.
	Given $p \in \dom{\deltaCOrd}$, if $p(0) = 0$
	we let $F(p) = 1\constantseq{0}$, i.e., a code
	for the tree $\extset\emptyseq$.
	If $p = 1q$, we let $F(p)$ be a code for the
	tree $T:= \extset\emptyseq \cup
	\intset{\seq{0}\concat\sigma}{\sigma \in
	\delta_\UT
	F(p)}$.
	Finally, if $p = 2q_0q_1\ldots$, we let $F(p)$
	be a code for the mix of the trees coded by the
	$F(q_n)$.
	It is now routine to check that $F \realizes
	\TreeWithRank$.
\end{proof}

\begin{cor}
	\th\label{witnessgeq}
	The operation
	$\WitnessgeqName:\mf{\LT\times\COrd\times\nat}{\UT}$,
	given by $S \in \Witnessgeq{\tree,\alpha,n}$
	iff
	\[
	\iPD{S}{\alpha} =
	\begin{cases}
		\extset{\emptyseq}, &\text{if } \enumTrees(n)
		\subseteq \iPD{\tree}{\alpha} \\
		\varnothing, &\text{otherwise,}
	\end{cases}
	\]
	is computable.
\end{cor}
\begin{proof}
	Given a labeled tree $\tree = (T,\varphi)$, a
	countable ordinal $\alpha$, and a natural
	number $n$, we output a tree $S$ of rank
	$\beta:= \min(\extset{\omega\cdot\alpha} \cup
	\intset{\rank_{T}(\sigma)}{\sigma \in
	\enumTrees(n)})$.

	We have $\iPD{S}{\alpha} = \extset\emptyseq$
	iff $\beta = \omega\cdot\alpha$ iff
	$\rank_{T}(\sigma) \geq \omega\cdot\alpha$ for
	each
	$\sigma \in \enumTrees(n)$ iff $\sigma \in
	\iPD{\tree}{\alpha}$ for each $\sigma \in
	\enumTrees(n)$ iff $\enumTrees(n) \subseteq
	\iPD{\tree}{\alpha}$, and $\iPD{S}{\alpha} =
	\varnothing$ otherwise.
\end{proof}

\begin{defi}
	We define $\graftName: \function{\LT \times \UT
	\times \UT \times \finbaire}{\LT}$ by
	\begin{enumerate}
		\item
			$\graft{\tree,S,U,\sigma} {\smallsetminus}
			\extensionSet{\finbaire}{\sigma} = \tree
			{\smallsetminus}
			\extensionSet{\finbaire}{\sigma}$

		\item
			$\concatenationSubtree{\graft{\tree,S,U,\sigma}}{\sigma}
			=
			\binaryMix{
			(\binaryTreeProduct{\concatenationSubtree{\tree}{\sigma}}{S})
			}{
			U
			}$
	\end{enumerate}
\end{defi}

\begin{defi}
	We define $\AuxName: \mf{\LT \times \LT \times
	\UT \times \COrd \times \finbaire \times
	\finbaire}{\LT}$ as follows.
	Given $\tree,\tree_\auxName \in \LT$, $\alpha
	\in \COrd$, and $\sigma,\tau \in \finbaire$
	such that $\length\sigma = \length\tau > 0$,
	let $\tree' \in
	\Aux{\tree,\tree_\auxName,U,\alpha,\sigma,\tau}$
	iff $\tree' =
	\graft{\tree,S_{\WitnessleqMark},\binaryTreeProduct{S_{\WitnessgeqMark}}{U},\sigma}$
	for some
	\[
	\begin{array}{rcl}
		S_{\WitnessleqMark}
		&\in&
		\Witnessleq{\tree_\auxName,\alpha,\last\tau}
		\\
		S_{\WitnessgeqMark}
		&\in&
		[\binaryTreeProduct{
		\finiteTreeProduct{n <
		\length\tau-1}{(\Witnessleq{\tree_\auxName,\alpha,\tau(n)})}
		]}{
		\Witnessgeq{\tree_\auxName,\alpha,\last\tau}
		}.
	\end{array}
	\]
\end{defi}

Recall from \th\ref{strongly_c_tightening} that
every computable or continuous multi-valued
function between represented spaces is
tightened by a strongly computable or strongly
continuous, respectively, multi-valued function
between the same spaces.
Therefore, in order to conclude that
$\iPD{\argument}{\alpha}$ is transparent for each
$\alpha$, it is enough to prove the following
stronger result.

\begin{thm}
	\th\label{derivative_uniformly_transparent}
	There is a computable operation
	$\transparentName: \mf{\cmf{\AT}{\AT} \times
	\COrd }{\cmf{\AT}{\AT}}$ such that $g \in
	\transparent{f}{\alpha}$ iff $\dom{f\iPD{
	\argument }{\alpha}} \subseteq \dom{g}$ and
	$\iPD{\mathcal{A}_g}{\alpha} \in
	f(\iPD{\mathcal{A}}{\alpha})$ for any
	$\mathcal{A} \in \dom{f\iPD{ \argument
	}{\alpha}}$ and $\mathcal{A}_g \in
	g(\mathcal{A})$.
\end{thm}
\begin{proof}
	Let $f \in \cmf{\AT}{\AT}$ be given in the form
	of a Turing machine $M$ which strongly computes
	$f$ with some given oracle $q$.
	Let $F:\pf\baire\baire$ be defined with
	$\dom{F} = \dom{f\delta_\AT}$ by letting $F(p)$
	be the output of $M$ on input
	$\tuple{p,\constantseq{0}}$ and oracle $q$.
	Thus $F$ is a computable realizer of $f$, so by
	\th\ref{positive_information} we can assume
	that for each $m$ there exists a
	computable subset $X_m \subseteq \nat$ such
	that $\enumTrees(m) \subseteq \delta_{\LT}
	F(p)$ iff $\enumTrees(n) \subseteq
	\delta_{\LT}(p)$ for some $n \in X_m$.
	Thus we can construct a computable labeled tree
	$\tree_F$ which \emph{represents} $F$, as
	follows.
	The nodes of length $1$ of $\tree_F$ are
	bijectively associated to the pairs $(n,\ell)$
	such that $\enumTrees(\ell)$ is a linear tree
	of height 1 and $n \in X_\ell$.
	If $\sigma \in \tree_F$ is associated to
	$(n,\ell)$, then
	\begin{enumerate}
		\item
			the label of $\sigma$ in $\tree_F$ is the
			label of the node of $\enumTrees(\ell)$ at
			height $\length\sigma$, and

		\item
			the children of $\sigma$ in $\tree_F$ are
			bijectively associated to the pairs
			$(n',\ell')$ such that $\enumTrees(\ell')$
			is a
			linear tree of height $\length\sigma+1$,
			$n' \in X_{\ell'}$, and $\enumTrees(\ell)
			\subseteq \enumTrees(\ell')$.
	\end{enumerate}

	It is now straightforward to check that if
	$\delta_\LT F(p)$ is not empty then it is
	bisimilar to the subtree $\tree_p$ of
	$\tree_F$ composed of the root plus those
	$\sigma$ which are associated to $(n,\ell)$
	with $\enumTrees(n) \subseteq \delta_\LT(p)$.

	Formally, our goal now is to computably define
	a Turing machine $M'$ from $M$, $q$, and
	$\alpha$, such that the function $g:=
	g_{M',q}$ from the proof of
	\th\ref{strongly_c_tightening} has the desired
	properties.
	To simplify the presentation, we will define
	$g$ directly and leave the definition of $M'$
	implicit.
	Thus, we want to define a computable
	$g:\mpf\LT\LT$ such that
	for any $p \in \dom{F}$ and any $\tree' \in
	g(\delta_\LT(p))$, letting $\delta_\LT(p') =
	\iPD{\delta_\LT(p)}{\alpha}$, we have:
 \begin{enumerate}
	 \item
		 if $\delta_\LT F(p) \neq \varnothing$ then
		 $\iPD{\tree'}{\alpha} \bisim \tree_{p'}$;

	 \item
		 if $\delta_\LT F(p) = \varnothing$ then
		 $\iPD{\tree'}{\alpha} = \varnothing$.
 \end{enumerate}

 Again, since $F$ is computable, there exists a
 computable $z \in \baire$ such that $\delta_\LT
 F(p) = \varnothing$ iff
 $\enumTrees(z(n)) \subseteq \delta_\LT(p)$ holds
 for some $n \in \nat$.
 Given $p \in \dom{f \iPD{ \argument }{\alpha}
 \delta_\AT}$, let $\tree:= \delta_\LT(p)$ and $U
 \in
 \WitnessAbsence(\tree,\alpha,z)$.
 Therefore, if $\delta_\LT F(p) \neq \varnothing$
 then $\enumTrees(z(m)) \not\subseteq
 \delta_\LT(p)$ for all $m \in \nat$ and
 thus $\iPD{U}{\alpha}$ is pruned and nonempty,
 and if $\delta_\LT F(p) = \varnothing$ then
 $\enumTrees(z(m)) \subseteq
 \delta_\LT(p)$ for some $m \in \nat$ and thus
 $\iPD{U}{\alpha} = \varnothing$.
 Let $V \in \TreeWithRank(\omega\cdot\alpha)$, so
 that $\iPD{V}{\alpha} = \extset\emptyseq$.
 Let $\tree_0 =
 \binaryTreeProduct{(\binaryMix{\tree_F}{V})}{U}$,
 so that if $\delta_\LT F(p) \neq \varnothing$
 then
 $\iPD{\tree_0}{\alpha} \bisim
 \iPD{\tree_F}{\alpha}$, and if $\delta_\LT F(p)
 = \varnothing$ then $\iPD{\tree_0}{\alpha} =
 \varnothing$.
 We let any node in $\tree_0$ coming from
 $\tree_F$ be associated to the same pair
 $(n,\ell)$ as the corresponding node in
 $\tree_F$.

 Now suppose we are at stage $s>0$ of the
 construction, so that we have already built a
 tree $\tree_{s-1}$.
 Let $\sigma:= \bij(s)$, where
 $\bij:\function{\nat}{\finbaire}$ is any
 computable bijection such that $\bij(s)
 \subseteq \bij(s')$
 implies $s \leq s'$.
 If $\sigma \not \in \tree_{s-1}$ or $\sigma \in
 \tree_{s-1}$ but is not associated to any
 $(n,\ell)$, then let $\tree_{s} =
 \tree_{s-1}$.
 Otherwise suppose $\restr{\sigma}{(m+1)}$ is
 associated to some $(n_m,\ell_m)$ for each $m <
 \length\sigma$.
 Let $\ast\sigma:=
 \seq{n_0,\ldots,n_{\length\sigma-1}}$ and define
 $\tree_{s}:=
 \Aux{\tree_{s-1},\tree,U,\alpha,\sigma,\ast\sigma}$.
 Recall that in this case we have
 \[
 \concatenationSubtree{\tree_{s}}{\sigma}
 =
 \binaryMix{
 (\binaryTreeProduct{\concatenationSubtree{\tree_{s-1}}{\sigma}}{S_{\WitnessleqMark}})
 }{
 (\binaryTreeProduct{S_{\WitnessgeqMark}}{U})
 }
 \]
 for some $S_{\WitnessleqMark},
 S_{\WitnessgeqMark}$ as in the definition of
 $\AuxName$.
 Hence we let each descendant $\sigma\concat\tau$
 of $\sigma$ in $\tree_s$ in which $\tau$ comes
 from the
 $\binaryTreeProduct{\concatenationSubtree{\tree_{s-1}}{\sigma}}{S_{\WitnessleqMark}}$
 component of $\binaryMixName$ above be
 associated to the same $(n,\ell)$ as the
 corresponding node in $\tree_{s-1}$.

 We then define $\tree'$ by letting $\sigma \in
 \tree'$ iff $\sigma \in \tree_s$ for $s =
 \preimage{\bij}(\sigma)$, with the label for
 $\sigma$ being its label in $\tree_s$ in the
 positive case.

 \begin{clm}
	 Every node of $\iPD{\tree'}{\alpha}$ other
	 than the root is associated to some pair
	 $(n,\ell)$.
 \end{clm}

 Indeed, it is easy to see that this is true of
 $\tree_0$.
 Thus if a node $\xi \in \iPD{\tree'}{\alpha}$ is
 not associated to some such pair, this means
 that $\xi$ was added to $\tree'$ at some
 stage $s > 0$ of the construction.
 Let $\sigma = \bij(s)$.
 In this case we have $\tree_{s}:=
 \Aux{\tree_{s-1},\tree,U,\alpha,\sigma,\ast\sigma}$,
 i.e.,
 $\tree_{s} =
 \graft{\tree_{s-1},S_{\WitnessleqMark},\binaryTreeProduct{S_{\WitnessgeqMark}}{U},\sigma}$
 for some
 \[
 \begin{array}{rcl}
	 S_{\WitnessleqMark}
	 &\in&\Witnessleq{\tree,\alpha,\last{\ast\sigma}}
	 \text{ and }\\
	 S_{\WitnessgeqMark} &\in& [\binaryTreeProduct{
	 \finiteTreeProduct{m <
	 \length\sigma-1}{(\Witnessleq{\tree,\alpha,\ast\sigma(m)})}
	 ]}{
	 \Witnessgeq{\tree,\alpha,\last{\ast\sigma}}
	 }.
 \end{array}
 \]
 The fact that $\xi$ is not associated to any
 pair $(n,\ell)$ implies that $\xi =
 \sigma\concat\eta$ for some $\eta \neq
 \emptyseq$
 coming from
 $\binaryTreeProduct{S_{\WitnessgeqMark}}{U}$.
 By construction the subtree of $\xi$ in $\tree'$
 is the same as in $\tree_{s}$, since for any $s'
 > s$ such that $\bij(s') \supseteq
 \xi$ we have $\tree_{s'} = \tree_{s'-1}$, and
 for any $s' > s$ such that $\sigma' = \bij(s')
 \not\supseteq \xi$ we have $\tree_{s'}
 {\smallsetminus}
 \extensionSet{\finbaire}{\sigma'} = \tree_{s'-1}
 {\smallsetminus}
 \extensionSet{\finbaire}{\sigma'}$.
 Hence
 $\concatenationSubtree{\iPD{\tree'}{\alpha}}{\xi}
 \subseteq
 \concatenationSubtree{\iPD{\binaryTreeProduct{S_{\WitnessgeqMark}}{U}}{\alpha}}{\eta}
 = \varnothing$, i.e., $\xi \not \in
 \iPD{\tree'}{\alpha}$.

 \begin{clm}
	 If $\delta_\LT F(p) = \varnothing$ then
	 $\iPD{\tree'}{\alpha} = \varnothing$.
 \end{clm}

 Indeed, if $\delta_\LT F(p) = \varnothing$ then
 $\iPD{U}{\alpha} = \varnothing$.
 Hence $\iPD{\tree_0}{\alpha} = \varnothing$, and
 at each stage $s > 0$ we either keep $\tree_s =
 \tree_{s-1}$, or else $\tree_s$
 differs from $\tree_{s-1}$ only in that
 \[
 \concatenationSubtree{\tree_{s}}{\sigma}
 =
 \binaryMix{
 (\binaryTreeProduct{\concatenationSubtree{\tree_{s-1}}{\sigma}}{S_{\WitnessleqMark}})
 }{
 (\binaryTreeProduct{S_{\WitnessgeqMark}}{U})
 }
 \]
 for some $S_{\WitnessleqMark},
 S_{\WitnessgeqMark}$ as in the definition of
 $\Aux{\tree_{s-1},\tree,U,\alpha,\sigma,\ast\sigma}$,
 where $\sigma = \bij(s)$.
 But then we have that
 \[
 \begin{array}{l}
	 \concatenationSubtree{\iPD{\tree_s}{\alpha}}{\sigma}
	 \\
	 \hspace{6em}=
	 \binaryMix{
	 (\binaryTreeProduct{\concatenationSubtree{\iPD{\tree_{s-1}}{\alpha}}{\sigma}}{\iPD{S_{\WitnessleqMark}}{\alpha}})
	 }{
	 \iPD{\binaryTreeProduct{S_{\WitnessgeqMark}}{U}}{\alpha}
	 } \\
	 \hspace{6em}=
	 \binaryTreeProduct{\concatenationSubtree{\iPD{\tree_{s-1}}{\alpha}}{\sigma}}{\iPD{S_{\WitnessleqMark}}{\alpha}},
 \end{array}
 \]
 so assuming by induction that
 $\iPD{\tree_{s-1}}{\alpha} = \varnothing$ holds,
 it follows that
 $\concatenationSubtree{\iPD{\tree_s}{\alpha}}{\sigma}
 = \varnothing$ as well.
 But then $\iPD{\tree_s}{\alpha} = \varnothing$,
 as desired.
 Therefore we have $\iPD{\tree'}{\alpha} =
 \varnothing$.

 For the rest of the proof we assume that
 $\delta_\LT F(p) \neq \varnothing$, which
 implies that $\iPD{U}{\alpha}$ is a pruned
 and nonempty tree.
 Furthermore, since $\iPD{V}{\alpha} =
 \extset\emptyseq$, we have $\emptyseq \in
 \iPD{\tree'}{\alpha}$.

 \begin{clm}
	 \th\label{survives}
	 Suppose $\sigma \in \tree' {\smallsetminus}
	 \extset\emptyseq$.
	 Then $\sigma \in \iPD{\tree'}{\alpha}$ iff
	 $\enumTrees(\ast\sigma(m)) \subseteq
	 \iPD{\tree}{\alpha}$ for each $m <
	 \length\sigma$.
 \end{clm}

 Let $s = \preimage{\bij}(\sigma)$.
 Suppose $\enumTrees(\ast\sigma(m)) \not\subseteq
 \iPD{\tree}{\alpha}$ for some $m <
 \length\sigma$.
 Let $s' =
 \preimage{\bij}(\restr{\sigma}{(m+1)})$.
 Note that $\iPD{S_\WitnessleqMark}{\alpha} =
 \iPD{\binaryTreeProduct{S_\WitnessgeqMark}{U}}{\alpha}
 = \varnothing$ for any
 $S_\WitnessleqMark \in
 \Witnessleq{\tree,\alpha,\ast\sigma(m)}$ and
 $S_\WitnessgeqMark \in
 \Witnessgeq{\tree,\alpha,\ast\sigma(m)}$.
 Thus we also have
 $\concatenationSubtree{\iPD{\tree_s}{\alpha}}{\restr{\sigma}{(m+1)}}
 = \varnothing$.
 Put together, and also considering the preceding
 claim, the two last statements imply $\tau \not
 \in \iPD{\tree_{s''}}{\alpha}$ for
 any $\tau \supseteq \restr{\sigma}{(m+1)}$.
 Thus $\sigma \not \in \iPD{\tree'}{\alpha}$.
 Conversely, suppose $\enumTrees(\ast\sigma(m))
 \subseteq \iPD{\tree}{\alpha}$ for every $m <
 \length\sigma$.
 Then for any
 \[
 S_{\WitnessgeqMark} \in [\binaryTreeProduct{
 \finiteTreeProduct{m <
 \length\sigma-1}{(\Witnessleq{\tree,\alpha,\ast\sigma(m)})}
 ]}{
 \Witnessgeq{\tree,\alpha,\last{\ast\sigma}}
 }
 \]
 we have
 $\iPD{\binaryTreeProduct{S_{\WitnessgeqMark}}{U}}{\alpha}
 = \iPD{S_{\WitnessgeqMark}}{\alpha} =
 \extset\emptyseq$.
 In particular
 it follows that the descendants of $\sigma$ in
 $\tree_s$ which are not associated to any
 $(n,\ell)$ already guarantee that $\sigma \in
 \iPD{\tree'}{\alpha}$, as desired.

 Let $p'$ be such that $\delta_\LT(p') =
 \iPD{\tree}{\alpha}$.

 \begin{clm}
	 The trees $\iPD{\tree'}{\alpha}$ and
	 $\tree_{p'}$ are bisimilar.
 \end{clm}

 Define $B \subseteq \iPD{\tree'}{\alpha} \times
 \tree_{p'}$ by letting $\sigma \mathrel{B} \tau$
 iff $\sigma = \tau = \emptyseq$ or
 $\length\sigma = \length\tau$,
 $\restr{\sigma}{n} \mathrel{B} \restr{\tau}{n}$
 for each $n < \length\sigma$, and $\sigma$ and
 $\tau$
 are associated to the same pair $(n,\ell)$.
 In order to verify that $B$ is a bisimulation,
 the only nontrivial properties to check are
 (\ref{bac}) and (\ref{for}).
 So suppose $\sigma \mathrel{B} \tau$ and for
 (\ref{bac}) let $\tau'$ be a child of $\tau$ in
 $\tree_{p'}$.
 Then $\tau'$ is associated to some $(n',\ell')$
 where $\enumTrees(n') \subseteq \delta_\LT(p') =
 \iPD{\tree}{\alpha}$.
 But then by construction $\sigma$ has some child
 $\sigma'$ in $\tree'$ associated to
 $(n',\ell')$.
 By \th\ref{survives} we have $\sigma' \in
 \iPD{\tree'}{\alpha}$, and $\sigma' \mathrel{B}
 \tau'$ follows.
 Finally, for (\ref{for}) let $\sigma'$ be a
 child of $\sigma$ in $\iPD{\tree'}{\alpha}$.
 Again by \th\ref{survives} we get that $\sigma'$
 is associated to some $(n',\ell')$ such that
 $\enumTrees(n') \subseteq
 \iPD{\tree}{\alpha}$.
 But then $\tau$ must have a child $\tau'$ in
 $\tree_{p'}$ which is also associated to
 $(n',\ell')$, and therefore $\sigma' \mathrel{B}
 \tau'$.

 Our assumption that $\delta_\LT F(p) \neq
 \varnothing$ implies that both
 $\iPD{\tree'}{\alpha}$ and $\tree_{p'}$ are
 nonempty
 trees, and $B \neq \varnothing$.
 Hence we have $\iPD{\tree'}{\alpha} \bisim
 \tree_{p'}$ as desired.
\end{proof}

\begin{thm}
	\th\label{derivative_cylinder}
	The operation $\iPD{\argument}{\alpha}$ is a
	transparent cylinder.
\end{thm}
\begin{proof}
	Transparency follows directly from
	\th\ref{derivative_uniformly_transparent}.
	To see that $\iPD{\argument}{\alpha}$ is a
	cylinder, given a code $p$ of an abstract tree
	$\mathcal{A}$, let $\mathcal{A}_p$ be the
	abstract tree obtained from $\mathcal{A}$ by
	changing each of its labels $\ell$ to
	$\tuple{1,\ell}$ plus adding an infinite path
	with
	induced label
	$\seq{\tuple{0,p(n)}}_{n\in\nat}$.
	Then $\iPD{\mathcal{A}_p}{\alpha}$ is obtained
	from $\iPD{\mathcal{A}}{\alpha}$ by the same
	change of labels as above plus the
	addition of the same infinite path.
	Now both $p$ and $\iPD{\mathcal{A}}{\alpha}$
	can easily be reconstructed from
	$\iPD{\mathcal{A}_p}{\alpha}$ without needing
	direct
	access to $p$; in other words, $\id_\baire
	\times \iPD{\argument}{\alpha} \leqsW
	\iPD{\argument}{\alpha}$.
\end{proof}

Let $\AT_\lin$ be the subspace of $\AT$ composed
of the linear abstract trees, and let
$\AT^*_\lin$ be the subspace of $\AT_\lin$
composed of the linear abstract trees which have
a unique infinite induced label.
The spaces $\AT_\fb$ and $\AT^*_\fb$ are defined
analogously for finitely branching trees.
Note that $\AT^*\lin$ is composed exactly of the
nonempty pruned linear trees.
Let $\Linearize_\fb$ be the restriction of
$\PDName$ to $\AT^*_\fb$, and note that
$\Linearize_\fb:\function{\AT^*_\fb}{\AT^*_\lin}$.

\begin{lem}
	The operation $\Linearize_\fb$ is
	Weihrauch-equivalent to $\lim$.
\end{lem}
\begin{proof}
	($\lim \leqW \Linearize_\fb$)
	Given $p \in \dom{\lim}$, we can build an
	abstract finitely branching tree whose induced
	labels are exactly the sequences of the form
	$\restr{(p)_n}{n}$.
	Since $\lim(p)$ is well defined, this tree is
	in the domain of $\Linearize_\fb$; applying
	this map to this tree results in a linear
	tree with an infinite branch labeled $\lim(p)$.

	($\Linearize_\fb \leqW \lim$)
	Given a name $p$ of an abstract tree
	$\mathcal{A}$ in the domain of
	$\Linearize_\fb$, let $\tree = \delta_\LT(p)$
	be one of its
	representatives.
	Since $\tree$ is bisimilar to a finitely
	branching tree, for each $\sigma \in \tree$ by
	K\H{o}nig's lemma we have that $\sigma \in
	\PD{\tree}$ iff
	$\concatenationSubtree{\tree}{\sigma}$ has
	infinite height.
	Therefore deciding whether $\sigma \in
	\PD{\tree}$ holds can be done with a single use
	of $\lim$, and since $\lim$ is parallelizable,
	one application of $\lim$ suffices to decide
	this for all $\sigma \in \tree$ at once.
	With this information we can construct
	$\Linearize_\fb(\mathcal{A})$.
\end{proof}

\begin{thm}
	The operation $\Linearize_\fb$ is transparent.
\end{thm}
\begin{proof}
	The proof is a simplified version of the proof
	of \th\ref{derivative_uniformly_transparent}.

	Let $f:\mpf{\AT^*_\lin}{\AT^*_\lin}$ be
	computable.
	Then $f$ has a realizer $F$ such that for each
	$\tau \in \finbaire$ there exists a computable
	$X_\tau \subseteq \finbaire$ such that
	$\tau$ is an induced label of the tree
	$\delta_\AT F(p)$ iff $\xi$ is an induced label
	of $\delta_\AT(p)$ for some $\xi \in
	X_\tau$.
	Let $p$ be given and $\tree:= \delta_\LT(p)$.
	We can computably define a labeled tree
	$\tree_G$ with the following properties.
	The nodes at level $1$ of $\tree_G$ are
	bijectively associated to the pairs
	$(\xi,\tau)$ such that $\length\tau=1$ and $\xi
	\in
	X_\tau$ is an induced label of $\tree$.
	Recursively, if $\sigma \neq \emptyseq$ is in
	$\tree_G$ and is associated to a pair
	$(\xi,\tau)$, then we have:
	\begin{enumerate}
		\item
			The induced label of $\sigma$ in $\tree_G$
			is $\tau$.

		\item
			\label{fb_transp_height_condition}
			If some node of $\tree$ with induced label
			$\xi$ has rank at least $\length\tau+1$,
			then the children of $\sigma$ in $\tree_G$
			are
			bijectively associated to the pairs
			$(\xi',\tau')$ such that $\length{\tau'} =
			\length\tau + 1$, $\tau' \supset \tau$, and
			$\xi'
			\in X_{\tau'}$ is an induced label of
			$\tree$; otherwise $\sigma$ is a leaf of
			$\tree_G$.
	\end{enumerate}

	\begin{clm}
		The trees $\PD{\tree_G}$ and $\tree_{FH}:=
		\delta_\LT F H (p)$ are bisimilar.
	\end{clm}

	Let $H \realizes \PDName:\function\LT\LT$.
	To see that $\PD{\tree_G} \bisim \tree_{FH}$,
	let $\sigma \mathrel{B} \tau$ iff $\sigma =
	\tau = \emptyseq$ or $\sigma$ and $\tau$
	have the same induced labels in $\PD{\tree_G}$
	and $\tree_{FH}$, respectively.
	Now suppose $\sigma \mathrel{B} \tau$, and let
	$\sigma$ be associated to $(\xi_0,\tau_1)$.
	Let $\sigma'$ be a child of $\sigma$ in
	$\PD{\tree_G}$.
	It follows that $\sigma'$ is associated to some
	pair $(\xi_1,\tau_1)$ such that $\tau_0
	\subseteq \tau_1$.
	Since $\sigma'$ is in the pruning derivative of
	$\delta_\LT G (p)$, by condition
	\ref{fb_transp_height_condition} of the
 construction it follows that there are nodes of
 $\tree$ of arbitrary length whose labels extend
 $\xi_1$.
 Since $\tree$ is bisimilar to a finitely
 branching tree, this implies that some node
 $\nu$ of $\tree$ with induced label $\xi_1$ is
 the root of a subtree of $\tree$ of infinite
 height.
 Thus $\nu$ is in $\delta_\LT H (p)$, and since
 $\xi_1 \in X_{\tau_1}$ it follows that some node
 with induced label $\tau_1$ is
 in $\tree_{FH}$.
 Finally, since $\tree_{FH}$ is linear, it
 follows that $\tau$ has a child $\tau'$ with
 induced label $\tau_1$, and thus $\sigma'
 \mathrel{B} \tau'$.
 Conversely, let $\tau'$ be a child of $\tau$ in
 $\tree_{FH}$, and let $\tau_1$ be its induced
 label.
 Therefore, some $\xi_1 \in X_{\tau_1}$ is an
 induced label in $\PD{\tree}$, and thus some
 node $\nu$ of $\tree$ has $\tau_1$ as its
 induced label and is the root of a subtree of
 $\tree$ of infinite height.
 This implies that some child $\sigma'$ of
 $\sigma$ in $\tree_G$ is associated to
 $(\xi_1,\tau_1)$, and that such $\sigma'$ is
 also in
 $\PD{\tree_G}$.
 Therefore $\sigma' \mathrel{B} \tau'$.
 Finally, note that $\tree_{FH}$ contains an
 infinite path since
 $f:\mpf{\AT^*_\lin}{\AT^*_\lin}$, which implies
 that $\tree_G$ and
 $\PD{\tree_G}$ also contain an infinite path.
 Therefore $B \neq \varnothing$ and
 $\PD{\tree_G}$ is bisimilar to $\tree_{FH}$.

 \begin{clm}
	 The tree $\tree_G$ is bisimilar to a finitely
	 branching tree.
 \end{clm}

 By construction, nodes of $\tree_G$ which have
 the same induced label have bisimilar (indeed,
 isomorphic) subtrees.
 Thus if some node $\sigma$ of $\tree_G$ has
 infinitely many children $\sigma_n$ which are
 roots of non-bisimilar subtrees, then the
 labels of the $\sigma_n$ are pairwise distinct.
 Therefore the $\sigma_n$ must be associated to
 elements $(\xi_n,\tau_n)$ such that the $\tau_n$
 are pairwise $\subseteq$-incomparable.
 But $\tree$ is bisimilar to a finitely branching
 tree; thus in particular only finitely many
 different labels occur on each of its
 levels.
 This implies that $\lim_{n \in \nat}
 \length{\tau_n} = \infty$, and therefore
 arbitrarily long prefixes of the infinite
 induced label
 of $\tree$ occur among the prefixes of the
 $\tau_n$.
 But then we cannot have that all $\xi_n$ have
 the same length $\length\sigma +1$, a
 contradiction.
\end{proof}

\begin{lem}
	The operation $\Linearize_\fb$ is a cylinder.
\end{lem}
\begin{proof}
	Given a name $p$ of an abstract tree
	$\mathcal{A} \in \dom{\Linearize_\fb}$, let
	$\mathcal{A}_p$ be the tree obtained from
	$\mathcal{A}$ by changing the label $\ell$ of
	any node $\sigma \neq \emptyseq$ to
	$\tuple{\ell,p(\length\sigma)}$.
	Then $\Linearize_\fb(\mathcal{A}_p)$ is
	obtained from $\Linearize_\fb(\mathcal{A})$ via
	the same transformation, and since
	$\Linearize_\fb(\mathcal{A})$ has an infinite
	branch, it is easy to reconstruct both $p$ and
	$\Linearize_\fb(\mathcal{A})$ from
	$\Linearize_\fb(\mathcal{A}_p)$.
	In other words, $\id_\baire \times
	\Linearize_\fb \leqsW \Linearize_\fb$.
\end{proof}

\subsection{Games for functions of a fixed Baire
class}

For an ordinal $\alpha = \lambda + 2n$, with
$\lambda$ a limit ordinal and $n$ a natural
number, let $\Linearize^{\alpha}$ be the
corestriction of $\iPD{\argument}{\lambda + n}$
to
$\AT^*_\lin$, let $\Linearize^{\alpha}_\fb$ be
the corestriction of $\iPD{\argument}{\lambda +
n}$ to $\AT^*_\fb$, and finally
let $\Linearize^{\alpha+1} = \Linearize_\fb
\comp \Linearize^{\alpha}_\fb$.

\begin{cor}
	\th\label{alpha_tree_game_parametrized}
	\label{alpha_tree_game_parametrized_page}
	Let $\alpha < \omega_1$.
	We have that $\Linearize^{\alpha}$ is a
	transparent cylinder which is
	Weihrauch-complete for the Baire class $\alpha$
	functions.
	Therefore the
	$(\Linearize^{\alpha},\Label)$-Wadge game
	characterizes the Baire class $\alpha$
	functions.
\end{cor}
\begin{proof}
	Suppose $\alpha = \lambda + 2n$.
	We have that $\iPD{\argument}{\lambda+n}$ is a
	transparent cylinder which is
	Weihrauch-complete for the Baire class
	$\lambda+2n$
	functions, so to see that the same holds for
	$\Linearize^{\lambda+2n}$, by
	\th\ref{corestr_of_transparent_cylinder} it is
	enough to
	show that $\AT^*_\lin$ strongly encodes
	$\baire$.
	But any $F:\mpf\baire\baire$ is easily seen to
	be strongly Weihrauch-equivalent to the map
	$F':\mpf\baire{\AT^*_\lin}$ which assigns
	$x \in \dom{F}$ to any linear abstract tree
	whose unique infinite label is in $F(x)$.
	Now suppose $\alpha = \lambda + 2n + 1$.
	Since $\AT^*_\lin \subseteq \AT^*_\fb \subseteq
	\AT$, by \th\ref{corestr_to_larger_subspace}
	and the fact that
	$\Linearize^{\lambda+2n}$ is Weihrauch-complete
	for Baire class $\lambda + 2n$ it follows that
	$\Linearize^{\lambda+2n}_\fb$ also has
	this property.
	Now, since $\Linearize_\fb$ is a transparent
	cylinder which is Weihrauch-complete for the
	Baire class 1 functions, the result follows.
\end{proof}

In other words, for $\alpha = \lambda + 2n$ with
$\lambda$ a limit ordinal and $n$ a natural
number, the restriction of the tree game in which
the final tree built by player $\II$ must have
$(\lambda + n)$\textsuperscript{th} pruning
derivative bisimilar to a linear tree
characterizes the Baire class $\alpha$
functions, and for $\alpha = \lambda + 2n + 1$,
the restriction of the tree game in which the
final tree built by player $\II$ must have
$(\lambda + n)$\textsuperscript{th} pruning
derivative bisimilar to a finitely branching tree
characterizes the Baire class
$\alpha$ functions.

\providecommand{\url}[1]{\texttt{#1}}
\providecommand{\urlprefix}{URL }
\bibliographystyle{alpha}
\bibliography{refs}

\end{document}